\documentclass[final]{amsart}
\usepackage{microtype}
\usepackage[foot]{amsaddr}

\usepackage{amsmath}
\usepackage{amssymb}
\usepackage{graphicx}
\usepackage{color}
\usepackage{xcolor}
\definecolor{dblue}{rgb}{0,0,0.3}
\usepackage{hyperref}
\usepackage{enumitem}
\usepackage{mathtools}
\usepackage{enumitem}
\usepackage{leftidx}
\usepackage{isomath}
\usepackage{cite}

\usepackage[small]{complexity}

\usepackage{algorithm}
\usepackage{algpseudocode}
\usepackage{varwidth}

\usepackage{tikz}
\usetikzlibrary{positioning}
\usetikzlibrary{calc}

\newcommand{\remove}[1]{}

\newtheorem{theorem}{Theorem}[section]
\newtheorem{definition}{Definition}
\newtheorem{proposition}{Proposition}[section]
\newtheorem{lemma}{Lemma}[section]
\newtheorem{claim}{Claim}[theorem]
\newtheorem{corollary}{Corollary}[section]
\newtheorem{example}{Example}[section]
\newtheorem{remark}{Remark}
\usepackage[normalem]{ulem}

\newtheorem{innercustomthm}{Theorem}
\newenvironment{customthm}[1]
  {\renewcommand\theinnercustomthm{#1}\innercustomthm}
  {\endinnercustomthm}

\newcommand{\N}{\mathbb{N}}
\newcommand{\n}{\{1,\ldots,n\}}
\newcommand{\m}{\{1,\ldots,m\}}

\newcommand{\f}{\bar{f}}

\newcommand{\x}{\bar{x}}
\newcommand{\y}{\bar{y}}
\newcommand{\ph}{\bar{\phi}}
\newcommand{\ba}{\bar{a}}
\newcommand{\bb}{\bar{b}}
\newcommand{\bc}{\bar{c}}

\newcommand{\un}{uniform}

\newcommand{\eat}[1]{}

\newcommand{\rMod}{\textrm{Mod}}
\newcommand{\rS}{\rm{\Sigma}}
\newcommand{\rP}{\rm{\Pi}}
\newcommand{\rD}{\rm{\Delta}}
\newcommand{\pr}{\rm{pr}}

\newcommand{\imp}{\rm{Imp}}

\newcommand{\scc}{\rm{scc}}

\DeclareMathOperator{\maj}{maj}

\title{On the Computational Complexity of Non-dictatorial Aggregation}

\author[L. Kirousis]{Lefteris Kirousis$^1$}
\email{lkirousis@math.uoa.gr}

\author[P. G. Kolaitis]{Phokion G. Kolaitis$^2$}
\email{kolaitis@ucsc.edu}

\author[J. Livieratos]{John Livieratos$^1$}
\email{jlivier89@math.uoa.gr}

\address{\hspace{-0.45cm}$^1$Department of Mathematics, National and Kapodistrian University of Athens\\$^2$Computer Science and Engineering Department, UC Santa Cruz and IBM Research - Almaden}

\begin{document}
\begin{abstract}
 We investigate when non-dictatorial aggregation is possible from an algorithmic perspective, where non-dictatorial aggregation means that the votes cast by the members of a society can be aggregated in such a way that there is no single member of the society that always dictates the collective outcome.  We consider the setting in which the members of a society take a position on a fixed collection of issues, where for each issue several different alternatives are possible, but the combination of choices must belong to a given set  $X$ of allowable voting patterns. Such a set  $X$ is called a possibility domain if there is an aggregator that is non-dictatorial, operates separately on each issue, and returns values among those cast by the society on each issue.
We design a polynomial-time algorithm that decides, given a set  $X$ of voting patterns, whether or not
 $X$ is a possibility domain. Furthermore, if $X$ is a possibility domain, then the algorithm constructs in polynomial time a non-dictatorial aggregator for $X$. Furthermore, we show that the question of whether a Boolean domain $X$ is a possibility domain is in \textrm{NLOGSPACE}. We also design a polynomial-time algorithm that decides whether $X$ is a uniform possibility domain, that is,  whether $X$ admits an aggregator that is non-dictatorial even when restricted to any two positions for each issue. As in the case of possibility domains, the algorithm also constructs in polynomial time a uniform non-dictatorial aggregator, if one exists. Then, we turn our attention to the case where $X$ is given implicitly, either as the set of assignments satisfying a propositional formula, or as a set of consistent evaluations of a sequence of propositional formulas. In both cases, we provide bounds to the complexity of deciding if $X$ is a (uniform) possibility domain. Finally, we extend our results to four types of aggregators that have appeared in the literature: generalized dictatorships, whose outcome is always an element of their input, anonymous aggregators, whose outcome is not affected by permutations of their input, monotone, whose outcome does not change if more individuals agree with it and systematic, which aggregate every issue in the same way.
\end{abstract}

\maketitle

\section{Introduction} \label{sec:intro}
The study of vote aggregation has occupied a central place in social choice theory. A broad framework for carrying out this study is as follows. There is a fixed collection of issues on each of which every member of a society takes a \emph{position}, that is, for each issue, a member of the society can choose between a number of alternatives. However, not every combination of choices is allowed, which means that the vector of the choices made by a member of the society must belong to a given set  $X$ of allowable voting patterns, called \emph{feasible evaluations}. The goal is to investigate properties of \emph{aggregators}, which are functions   that take as  input the votes cast by the members of the society and return as output a feasible evaluation that represents the collective position of the society on each of the issues at hand. A concrete key problem studied in this framework is to determine whether or not a \emph{non-dictatorial} aggregator exists, i.e., whether or not it is possible to aggregate votes in such a way that no individual member of the society dictates the votes on the society. A set  $X$ of feasible evaluations is called a \emph{possibility domain} if it admits a non-dictatorial aggregator; otherwise, $X$ is called an \emph{impossibility domain}.
This framework is broad enough to account for several well-studied cases of vote aggregation, including the case of \emph{preference} aggregation for which J. K. Arrow \cite{arrow1951social} established his celebrated impossibility theorem  and the case of \emph{judgment} aggregation
%phk (see for example \cite{DBLP:reference/choice/Endriss16} introductory exposition).
(for a survey of the latter, see \cite{DBLP:reference/choice/Endriss16}).

The investigation of the existence of non-dictatorial aggregators is typically carried out under two assumptions: (a) the aggregators are \emph{independent} and (b) the aggregators are \emph{conservative} (also known as \emph{supportive} or \emph{grounded}). The independence assumption means that the aggregator is an issue-by-issue aggregator, so that an independent aggregator on $m$ issues can be identified with an $m$-tuple $(f_1,\ldots,f_m)$ of functions aggregating the votes on each issue. It is not difficult to see that, if independence is not assumed,  then non-dictatorial aggregation is trivially always possible. The notion of independence corresponds to that  of \emph{independence of irrelevant alternatives} (IIA)  used by J. K. Arrow \cite{arrow1951social} in the preference aggregation framework, under a natural embedding (provided by Dokow and Holzman in \cite{dokow2009aggregation})  of the preference aggregation framework to the framework  described here. The conservativeness (or supportiveness) assumption means that, for every issue, the position returned by the aggregator is one of the positions held by the members of the society on that issue. An immediate consequence of this property, is that the social decision cannot have ``compromises'' between the choices of the voters.

By now, there is a body of research on identifying criteria that characterize when a given set $X$ of feasible evaluations is a possibility domain, that is, it admits an $m$-tuple of $n$-ary conservative functions that are not all projections to the same coordinate. The first such criterion was established by Dokow and Holzman \cite{dokow2010aggregation} in the Boolean framework, where, for each issue, there are exactly two alternatives (say, $0$ and $1$) for the voters to choose from. Specifically, Dokow and Holzman showed that a set $X\subseteq \{0,1\}^m$ is a possibility domain if and only if $X$ is affine or $X$ is not totally blocked.  Informally, the notion of \emph{total blockedness}, introduced first by Nehring and Puppe \cite{nehring2002stategy},  asserts that the social position on any issue can be inferred from  the social position on any issue. This notion appears in many of the characterization results for possibility domains. In the non-Boolean framework (where, for some issues, there may be more than two alternatives), Dokow and Holzman \cite{dokow2010aggregationnonB} extended the notion of total blockedness and used it to give a sufficient condition for a set $X$  to be a possibility domain. Szegedy and Xu \cite{szegedy2015impossibility} used tools from universal algebra to characterize when a totally blocked set $X$ is a possibility domain. The reason for this is, as Kirousis et al. \cite{kirousis2019aggregation} showed, that total-blockedness corresponds to a weak notion of impossibility, specifically a set $X$ admits a binary non-dictatorial aggregator if and only if $X$ is not totally blocked. A consequence of these results is that a set  $X$  is a possibility domain if and only if $X$ admits a binary non-dictatorial aggregator or a ternary non-dictatorial aggregator; in other words, non-dictatorial aggregation is possible for a society of some size  if and only if it is possible for a society with just two members or with just three members. This line of work was pursued further by Kirousis et al. \cite{kirousis2019aggregation}, who characterized possibility domains in terms of the existence of binary non-dictatorial aggregators or ternary non-dictatorial aggregators of a particular form.

All the aforementioned results are situated in what Dokow and Holzman call the \emph{abstract} framework, where the domain $X$ is given explicitly as a set of $m$-ary tuples. Earlier works in this framework include, but are not limited to, \cite{fishburn1986aggregation,kasher1997question,rubinstein1986algebraic}. Judgment aggregation has been extensively studied in the framework of the \emph{logic-based} approach (see \cite{list2002aggregating,list2009judgment}). In that approach, there is a tuple $\ph=(\phi_1,\ldots,\phi_m)$ of propositional formulas, called the \emph{agenda}, and the set $X$ of feasible evaluations comprises consistent judgments concerning the validity of the formulas. A judgment over a formula $\phi$ amounts to deciding if $\phi$ or $\neg\phi$ is true and it is \emph{consistent} if of all the formulas deemed true can be satisfied simultaneously, by at least one assignment of values. Note that in this setting, the allowed positions on each issue are necessarily two and thus it is situated in the Boolean framework. In this framework, Endriss et al. \cite{endriss2012complexity} have studied the computational complexity of three interesting problems: the \emph{winner determination problem}, where we want to decide if a specific formula belongs to the collective decision on a given agenda under some aggregation procedure, the \emph{strategic manipulation} problem, where we want to decide if an agent has an incentive to misrepresent his or her preferences in order to change the collective decision on a given agenda under some aggregation procedure, and the \emph{safety of the agenda} problem, where we want to decide if a class of aggregators preserves the logical consistency restrictions of an agenda. This last problem is related to our framework, with the difference that we search for at least one aggregator of a given class preserving the logical restrictions of our domain. Also, \cite{terzopoulou2018aggregating} have studied the case where the individuals do not need to make a decision on every issue of the agenda, but can instead provide partial judgments.

A different approach is the one used by Grandi and Endriss \cite{grandi2013lifting}, where the set of restrictions is provided by a propositional formula $\phi$, called an \emph{integrity constraint}, as the set of its satisfying assignments. The agents here are required in a way to select, in a consistent manner, which variables of the formula are true and which are not, instead of doing so for entire propositional formulas of an agenda. Endriss et al. \cite{endriss2016succinctness} study the relation of this framework with the logic-based one, in terms of the \emph{succinctness} of the corresponding languages used and of the computational complexity of the winner determination problem.

\paragraph{Summary of Results.} The aforementioned investigations have characterized possibility domains (in both the Boolean  and the non-Boolean frameworks) in terms of \emph{structural} conditions.
Our goal is to investigate possibility domains using the \emph{algorithmic lens} and, in particular, to study the following algorithmic problem: given a set $X$ of feasible evaluations, determine whether or not $X$ is a possibility domain. Szegedy and Xu \cite[Theorem 36]{szegedy2015impossibility} give algorithms for this problem, but  these  algorithms have very high running time; in fact, they run in exponential time in the number of issues and in the number of positions over each issue, even when confined to the  Boolean framework.

We design a polynomial-time algorithm that, given a set $X$ of feasible evaluations (be it in the Boolean or the non-Boolean framework),  decides whether  $X$ is a possibility domain. Furthermore, if $X$ is a possibility domain, then the algorithm  produces a binary non-dictatorial  or a ternary non-dictatorial aggregator for $X$.

The first step towards this result is to show that there is a polynomial-time algorithm that given a set $X$ of feasible evaluations, decides whether  $X$ admits a binary non-dictatorial aggregator  (as mentioned earlier, this amounts to $X$ not being totally blocked).  In fact, we show a stronger result, namely, that this problem
is expressible in \emph{Transitive Closure Logic (TCL)}, an extension of \emph{first-order logic}  with the \emph{transitive closure} operator; see \cite{DBLP:books/sp/Libkin04} for the precise definitions. As a consequence, the problem of deciding whether $X$ admits a binary non-dictatorial aggregator
 is in \textrm{NLOGSPACE}. Using this result, we then show that the problem of deciding whether a set $X\subseteq\{0,1\}^m$ is a possibility domain is  in $\textrm{NLOGSPACE}$.

 After this, we turn our attention to \emph{{\un} possibility domains}, which were introduced by Kirousis et al. \cite{kirousis2019aggregation} and which form a proper subclass of the class of possibility domains. Intuitively, {\un} possibility domains are sets of feasible evaluations that admit an aggregator that is non-dictatorial even when restricted to any two positions for each issue. This requirement forces the aggregating procedure to treat each Boolean (yes/no) decision of every issue in a non-dictatorial way. In that way, we can avoid various trivial cases of non-dictatorial aggregation, where different ``dictators'' are chosen for different sets of positions of an issue. Kirousis et al. established a tight connection between uniform possibility domains and a variant of constraint satisfaction problems by showing that  multi-sorted conservative constraint satisfaction problems are tractable on uniform possibility domains, whereas such constraint satisfaction problems are NP-complete on all other domains. Here, using the result by Carbonnel \cite{carbonnel2016dichotomy}, we give a polynomial-time algorithm for the following decision problem: given a set $X$ of feasible evaluations (be it in the Boolean or the non-Boolean framework), determine whether or not $X$ is a uniform possibility domain; moreover, if $X$ is a uniform possibility domain, then the algorithm produces a suitable uniform non-dictatorial aggregator for $X$.

We also study the problems of non-dictatorial and uniform non-dictatorial aggregation in case $X$ is Boolean and provided via an integrity constraint or by an agenda. In both cases, we provide bounds for the computational complexity of deciding if $X$ is a (uniform) possibility domain. Finally, we extend our results to three types of aggregators that have been used in the literature: generalized dictatorships \cite{grandi2013lifting}, \cite{diaz2019syntactic}, whose outcome is always an element of their input, anonymous aggregators, whose outcome is not affected by permutations of their input and monotone aggregators, whose outcome does not change if more individuals agree with it (see \cite{nehring2010abstract}) and, finally, systematic aggregators, which aggregate every issue in the same way. This last type of aggregators corresponds to \emph{polymorphisms}, that is functions that preserve domains, which have been extensively studied in the literature (e.g., in \cite{szendrei1986clones}) and which have been successfully used in matters of computational complexity (e.g., in \cite{bohler2003playing,bohler2004playing}) and aggregation theory (e.g., in \cite{szegedy2015impossibility,kirousis2019aggregation}).

The results reported here contribute to the developing field of computational social choice and pave the way for further exploration of algorithmic aspects of vote aggregation. In a sense, the question we investigate is the following: given a specific voting scheme, where some pre-defined rules of logical consistency apply, how difficult is it to decide if it is possible to design an aggregation rule with some desired properties, and then construct it in case it is? It should be noted that in the field of Judgment Aggregation, quite frequently the domains that impose the logical consistency restrictions are fixed. In such a setting, the algorithmic approach does not have much to offer, since in this case we have a one-off problem. On the other hand, it is not difficult to imagine groups of people that constantly need to make collective decisions over different sets of issues, where the logical restrictions that apply change according to the dependencies between the issues. In this scenario, algorithms that can quickly decide which aggregation rules can be applied could be useful.

In terms of the applicability of our algorithms in the abstract setting, there is an issue with the size of the input, since a domain $X$ given explicitly as a set of $m$-ary vectors, can be too large for practical purposes. However, as small or large the domain might be, the search space of its possible aggregators is exponentially larger, even provided characterization results like the ones in \cite{dokow2009aggregation,kirousis2019aggregation,szegedy2015impossibility} that restrict the search to binary or ternary aggregators. Thus, our tractability results in the abstract framework should be interpreted as showing that, given access to the domain, one encounters no more problems in deciding whether non-dictatorial aggregation is possible. Furthermore, the algorithm for finding and producing binary non-dictatorial aggregators (or deciding their lack thereof), is used in obtaining parts of the complexity upper bounds in the cases where the domain is given in compact form.
\vspace{0.3cm}

\textbf{Structure of the paper.} Theorems, lemmas, and corollaries are numbered per section, each referring to its own type, whereas claims are numbered within theorems.\begin{itemize}
    \item In Section \ref{sec:prelim}, we formally describe the abstract framework and the necessary tools
%phk in order
used to obtain the polynomial-time algorithms for deciding if a domain $X$ is a (uniform) possibility domain. Specifically, in Subsection \ref{sub:possibility}, we state the main theorems we use to obtain our results concerning possibility domains, i.e., the characterizations of possibility domains in the Boolean framework [Theorem \ref{thm:boolean}], totally blocked domains [Theorem \ref{thm:tot-block}], and possibility domains in the non-Boolean case [Theorem \ref{thm:pd}]; we also  prove Lemma \ref{lem:cartprod} that outlines the type of aggregators for domains that are Cartesian products. In Subsection \ref{sub:upd}, we define uniform possibility domains and state Theorem \ref{thm:updcar} that characterizes them. In Subsection \ref{sub:PH}, we discuss the \emph{Polynomial Hierarchy}.
\item In Section \ref{sec:abstract}, we present our polynomial time algorithms. Specifically, in Subsection \ref{ssec:poss-tract}, we present polynomial time algorithms that check whether  a domain admits binary non-dictatorial aggregators [Theorem \ref{thm:binary-tract}],  whether   a domain is totally blocked [Corollary \ref{cor:tot-block}], and whether a domain is  a possibility domain [Theorem \ref{thm:pdmain}; first main result]. We also show that all domains admit \emph{maximum symmetric aggregators} [Lemma \ref{lem:maxsym}] and give a polynomial time algorithm that produces one [Corollary \ref{cor:maxsym}]. In Subsection \ref{ssec:tcs}, we describe \emph{Transitive Closure Logic} and \emph{Least Fixed Point Logic}, showing that detecting whether a Boolean domain admits binary non-dictatorial aggregators is in TCL and in LFP [Theorem \ref{binary-TCL:cor}], that checking whether a Boolean domain is affine is in {\rm LOGSPACE} [Lemma \ref{lem:affinelogspce}], and that checking whether it is a possibility domain is in {\rm NLOGSPACE} [Theorem \ref{thm:boolean-tract}; second main result]. In Subsection \ref{ssec:upd-tract}, we present a polynomial time algorithm that checks whether a domain is a uniform possibility domain [Theorem \ref{thm:tractupd}; third main result], using a polynomial time algorithm by Carbonnel \cite{carbonnel2016dichotomy} that checks whether a domain admits polymorphisms with specific properties [Theorem \ref{thm:carb-tract}]. We also state an alternative characterization of local possibility domains in the Boolean framework [Corollary \ref{cor:lpdchar}] and provide a polynomial time algorithm that checks whether a domain is a local possibility domain [Corollary \ref{cor:tractlpd}].
\item In Section \ref{sec:implicit}, we present some complexity results for the case where $X$ is given implicitly, either by an integrity constraint or an agenda. Specifically, in Subsection \ref{ssec:logicbased}, we describe the logic-based framework, defining integrity constraints and agendas. In Subsection \ref{ssec:integr}, we define the domain $\mathrm{Imp}$,  show that it is an impossibility domain [Lemma \ref{lem:imp}], and prove that the problem of deciding if the domain of an integrity constraint admits non-dictatorial aggregators is in $\rS_2^P\cap\rP_2^P$ and $\coNP$-hard [Theorem \ref{thm:possintcon}; fourth main result], whereas the same problem for locally non-dictatorial aggregators is in $\rS_2^P\cap\rP_2^P$ and $\coNP$-hard [Theorem \ref{thm:locposintcon}; fifth main result]. In Subsection \ref{ssec:agenda}, we show that the problem of deciding whether the domain of an agenda admits non-dictatorial aggregators is in $\rD_3^P$ and $\coNP$-hard [Theorem \ref{pdagenda}; sixth main result], whereas the same problem for locally non-dictatorial aggregators is in $\rD_3^P$ and $\coNP$-hard [Theorem \ref{agendalocal}; seventh main result]. In Subsection \ref{ssec:other}, we state a characterization of Boolean domains admitting aggregators that are not generalized dictatorships [Theorem \ref{thm:gendict}] and prove that the problem of deciding whether the domain of an integrity constraint (respectively, agenda) admits aggregators that are not generalized dictatorships is in $\rS_2^P\cap\rP_2^P$ (respectively, $\rD_3^P$) and $\coNP$-hard [Corollaries \ref{cor:gendictup} and \ref{cor:gendictagendaup}]. We provide the exact same results for anonymous [Corollaries \ref{cor:anon} and \ref{results:anon}] and monotone aggregators [Theorem \ref{thm:monotone} and Corollary \ref{results:monotone}]. Finally, we show that deciding whether the domain of an integrity constraint admits systematic aggregators is $\coNP$-complete [Proposition \ref{icshaefer}], whereas that of an agenda is in $\rP_2^P$ and $\coNP$-hard [Proposition \ref{agendashaefer}].
\end{itemize}
\vspace{0.3cm}

\textbf{Relation to the conference version.} A preliminary version of this paper appeared in the Proceedings of the 17th International Conference on Relational and Algebraic Methods in Computer Science (RAMiCS 2018) \cite{DBLP:conf/RelMiCS/KirousisKL18}. The  present full version, in addition to detailed proofs and several improvements in the presentation, contains the results about expressibility in Transitive Closure Logic and membership in \textrm{NLOGSPACE}. Furthermore, all the results in Section \ref{sec:implicit} concerning the case in which $X$ is given implicitly are new.

\section{Preliminaries and Earlier Work}\label{sec:prelim}
In Subsection \ref{sub:possibility} we consider possibility domains both in the Boolean and non-Boolean case, whereas in Subsection \ref{sub:upd} we turn our attention to uniform possibility domains.

\subsection{Possibility Domains}\label{sub:possibility}
 Let $I=\{1,\ldots,m\}$ be a  set of issues. Assume that  the \emph{possible position values} of an individual (member of a society) for each issue are given by the finite set $A$. We also assume that $A$ has cardinality at least $2$. If $|A|=2$, we say that we are in the \emph{Boolean} framework; otherwise we say that we are in the \emph{non-Boolean} framework.

An \emph{evaluation} is an element of $A^m$. Let $X\subseteq A^m$ be a set of \emph{permissible} or \emph{feasible} evaluations.
 %For technical reasons,
To avoid degenerate cases,  we assume that for each $j=1,\ldots,m$, the $j$-th projection $X_j$ of $X$ has cardinality \emph{at least} $2$. Notationally, we refer to the elements of $X$ (or sometimes of $A^m$) as vectors, patterns, evaluations or assignments, depending on the context.

Let $n\geq 2$ be the number of individuals. We view each element $x$ of $X^n$  as  a $n\times m$ matrix that represent the choices of all individuals over every issue. The element $x_j^i$ of such a matrix $x$ will be the choice of the $i$-th individual over the $j$-th issue, for $i=1,\ldots,n$ and  $j=1,\ldots,m$. The $i$-th row $x^i$ will represent the choices of the $i$-th individual over every issue, $i=1,\ldots,n$,  and the $j$-th column $x_j$ the choices of every individual over the $j$-th issue, $j=1,\ldots,m$.

In order to \emph{aggregate} a set of $n$ feasible evaluations, we use $m$-tuples $F=(f_1,\ldots,f_m)$ of $n$-ary functions, where $f_j:X_j^n\to X_j$, $j=1,\ldots,m$. Such a $m$-tuple $F$  of functions will be called an ($n$-ary) \emph{aggregator} for $X$ if it satisfies:
\begin{enumerate}
\item \emph{Collective rationality:} $(f_1(x_1),\ldots,f_m(x_m))\in X$.
\item \emph{Conservativity:} $f_j(x_j)\in \{x_j^1,\ldots,x_j^n\}$, for all $j\in \m$.
\end{enumerate}
An aggregator  $F = (f_1, \ldots, f_m)$ is called  {\em dictatorial on $X$} if it always outputs a specific vector of its input, that is, if there is a number $d \in \{1, \ldots, n\}$ such that $(f_1,\ldots, f_m) = ({\rm pr}_d^n,\ldots,{\rm pr}_d^n)$, where ${\rm pr}_d^n$ is  the $n$-ary projection on the $d$-th coordinate; otherwise, $F$ is called  {\em non-dictatorial on $X$}. We say that $X$ {\em has a non-dictatorial aggregator} if, for some $n\geq 2$, there is  a $n$-ary non-dictatorial  aggregator on $X$. In such a case, we say that $X$ is a {\em possibility domain}. Otherwise, it is  an {\em impossibility domain}. A possibility domain is, by definition,  one where non-dictatorial aggregation is possible for societies of some cardinality, namely, the arity of a non-dictatorial aggregator.

Observe that if $X$ has a non-dictatorial aggregator of arity $n$, then it has also non-dictatorial aggregators of arity $k$, for all $k\geq n$. Indeed, given an $n$-ary non-dictatorial aggregator $F=(f_1,\ldots,f_m)$ for $X$, we can define, for any $k\geq n$, the aggregator $G=(g_1,\ldots,g_m)$ whose components $g_j$ simply ignore the last $k-n$ elements of their input and behave like the corresponding $f_j$'s. On the other hand, Dokow and Holzman \cite{dokow2010aggregationnonB} provide an example of a domain admitting a ternary non-dictatorial aggregator, but no binary ones.

The notion of an aggregator is akin to, but different from,  the notion of a polymorphism -- a fundamental notion in universal algebra; see, e.g., \cite{szendrei1986clones}. A polymorphism is in fact a systematic aggregator, that is, an aggregator all of whose components are the same operation. Let $A$ be a finite non-empty set. A \emph{constraint language} over $A$  is a finite set $\Gamma$ of relations of finite arities. Let $R$ be an $m$-ary relation on $A$. A  function $f:A^n\to A$ is a \emph{polymorphism} of $R$ if the following condition holds:
$$\text{if } x^1,\ldots,x^n\in R, \text{ then } (f(x_1),\ldots,f(x_m))\in R,$$
where $x^i=(x_1^i,\ldots,x_m^i)\in R$, $i=1,\ldots,n$ and $x_j=(x_j^1,\ldots,x_j^n)$, $j=1,\ldots,m$.
In this case, we also say that $R$ \emph{is closed under} $f$ or that $f$ \emph{preserves} $R$.
Finally, we say that $f$ is a \emph{polymorphism of a constraint language} $\Gamma$ if $f$ preserves every relation $R\in \Gamma$.

In fact, polymorphisms of a domain $X$ correspond to \emph{systematic} aggregators for $X$, that is, aggregators $F=(f_1,\ldots,f_m)$ where $f_j=f$, $j=1,\ldots,m$. In what follows, we denote $m$-tuples comprised of the same function $f$ by $\f$.

A function $f:A^n\to A$ is \emph{conservative} if, for all $a_1,\ldots,a_n\in A$, we have that  $f(a_1,\ldots,a_n)\in \{a_1,\ldots,a_n\}.$
Clearly, if $f:A^n\to A$ is a conservative polymorphism of an $m$-ary relation $R$ on $A$, then
the $m$-tuple $\bar{f} = (f, \ldots, f)$ is an $n$-ary aggregator for $R$.

We say that a ternary operation $f: A^3 \to A $ on an arbitrary set $A$ is  a {\em majority} operation if for all $x$ and $y$ in $A$,
$$f(x,x,y) = f(x,y,x) = f(y,x,x) =x;$$ we  say that  $f$ is a {\em minority} operation  if  for all $x$ and $y$ in $A$,
$$f(x,x,y) = f(x,y,x) = f(y,x,x) =y.$$

We  also say that a set $X$ of feasible evaluations \emph{admits a majority (respectively, minority) aggregator} if it admits a ternary ($n=3$) aggregator every component of which is a majority (respectively, minority) operation. Clearly,
$X$ admits a majority aggregator if and only if  there is a ternary aggregator $F = (f_1, \ldots, f_m)$ for $X$ such that, for  all $j=1, \ldots, m$
and for  all two-element subsets $B_j \subseteq X_j$, we have that
$f_j{\restriction{B_j}} = {\rm maj},$ where
\begin{equation*}
	{\rm maj} (x, y, z) = \begin{cases}
 x & \text{ if }  x= y \text{ or } x=z,\\
 y & \text{ if } y = z.
 \end{cases}
\end{equation*}
Similarly, $X$ admits a minority aggregator if and only if there is a ternary aggregator $F = (f_1, \ldots, f_m)$ for $X$ such that, for  all $j=1, \ldots, m$
and for  all two-element subsets $B_j \subseteq X_j$, we have that $f_j{\restriction{B_j}} = \oplus,$ where
\begin{equation*}
	\oplus (x, y, z) = \begin{cases}
 z & \text{ if }  x= y,\\
 x & \text{ if } y = z,\\
 y & \text{ if } x = z.
 \end{cases}
\end{equation*}

In the Boolean framework, a set  $X\subseteq \{0,1\}^m$  admits a majority aggregator if and only if the majority operation ${\rm maj}$ on $\{0,1\}^3$ is a polymorphism of $X$.
\begin{lemma}{\cite[Lemma 3.1B]{schaefer1978complexity}}\label{sch:bijunctive}
The majority operation ${\rm maj}$ on $\{0,1\}^3$ is a polymorphism of a set $X\subseteq\{0,1\}^m$ if and only if $X$ is a bijunctive logical relation,  i.e., $X$  is the set of satisfying assignments of a 2CNF-formula.
\end{lemma}
A set $X\subseteq \{0,1\}^m$ admits a minority aggregator if and only if the minority operation $\oplus$ on $\{0,1\}^3$ is a polymorphism of $X$.
\begin{lemma}{\cite[Lemma 3.1A]{schaefer1978complexity}}\label{sch:affine}
The minority operation $\oplus$ on $\{0,1\}^3$ is a polymorphism of a set $X\subseteq\{0,1\}^m$ if and only if $X$ is an affine logical relation, i.e., $X$  is the set of solutions of a system of linear equations over the two-element field.
\end{lemma}
Both these results are part of the proof in Schaefer's Dichotomy Theorem for the Satisfiability Problem, in \cite[Theorem 2.1]{schaefer1978complexity}. Since aggregators are by definition conservative, we identify from now on the ternary majority and minority aggregators with $\overline{\maj}$ and $\bar{\oplus}$ respectively.
  \begin{example}\label{ex:bijaff} Consider the sets $X_1$ and $X_2$ below.
  \begin{itemize}
      \item[(i)] The set $X_1 =\{0,1\}^3\setminus \{(1,0,1), (0,0,1), (0,0,0)\}$ is bijunctive, since it is the set of satisfying assignments of the 2CNF-formula $(x\lor y) \land (y \lor \neg z)$.
      \item[(ii)] The set $X_2=\{(0,0,1), (0,1,0), (1,0,0), (1,1,1)\}$ is affine, since it is the set of solutions of the equation $x+y+z=1$ over the two-element field.
      \item[(iii)] Both sets $X_1$ and $X_2$ are possibility domains, since $X_1$ admits a majority aggregator and $X_2$ admits a minority aggregator.
  \end{itemize}
\end{example}
Finally, a binary ($n=2$) aggregator is called a \emph{projection} aggregator if all its components are projections and symmetric, if all its components are symmetric, i.e., $f_j(x,y)=f_j(y,x)$, for all $x,y\in X_j$ and  for $j=1,\ldots,m$.

We now present two theorems that characterize possibility domains in the Boolean and the non-Boolean framework.
They   are the stepping stones towards
showing that the following decision problem is solvable in polynomial time:
 given a set $X$ of feasible evaluations, is $X$ a possibility domain?

\begin{customthm}{A}{\cite[Theorem 2.2]{dokow2010aggregation}}\label{thm:boolean}
Let $X\subseteq \{0,1\}^m$ be a set of feasible evaluations. The following two statements are equivalent.
\begin{enumerate}
\item $X$ is a possibility domain.
\item $X$ is affine or $X$ is not totally blocked.
\end{enumerate}
\end{customthm}

As mentioned in the Introduction, the existence of a binary non-dictatorial aggregator on $X$ is closely related to the total blockedness of $X$. We follow closely the notation and terminology used in \cite{dokow2010aggregationnonB}.

Let $X$ be a set of feasible voting patterns. Given subsets  $B_j \subseteq X_j, j=1, \ldots, m$, a \emph{sub box} is the product $B = \prod_{j=1}^m B_j $. It is called a $2$- sub-box if $|B_j| =2$, for all $j$. Elements of a box $B$ that belong also to $X$  will  be called {\em feasible evaluations within $B$}.

For a subset $K\subseteq\{1,\ldots,m\}$ and a tuple $x\in\prod_{j \in K}B_j$ is a {\em   feasible partial evaluation within $B$} if there exists a feasible evaluation $y$ within $B$ that extends $x$, i.e., $x_j = y_j$, for all $j \in K$. Otherwise,  we say that $x$  is  an \emph{infeasible} partial evaluation within $B$. $x$ is a {\em $B$-Minimal Infeasible  Partial Evaluation} ($B$-MIPE) if (i) it is an  infeasible  partial evaluation within $B$ and (ii) if for every $j \in K$,  there is a $b_j \in B_j$ such  that changing the  $j$-th  coordinate of $x$ to $b_j$ results into a feasible  partial evaluation within $B$.

We can now define a directed graph $G_X$ as follows. The  vertices of $G_X$ are the pairs of {\em distinct} elements $u,u'$ in $X_j$, for all $j =1, \ldots m$, denoted by $uu'_j$. Two vertices $uu'_k, vv'_l$ with $k \neq l$ are connected by a directed edge from  $uu'_k$ to $vv'_l$ if there exists a 2-sub-box $B = \prod_{j=1}^m B_j $, a set $K \subseteq \{1, \ldots, m\}$, and a $B$-MIPE $x = (x_j)_{j\in K}$ such that $k,l \in K$ and $B_k = \{u,u'\}$ and $ B_l =\{v,v'\}$ and $x_k =u$ and $x_l = v'$.

\begin{definition}{\cite{dokow2010aggregationnonB}}\label{def:tot-blockr}
We say that $X$ is {\em totally blocked} if the graph $G_X$ is strongly connected, i.e.,  every two distinct vertices $uu'_k, vv'_l$  are connected by a directed path (this must hold even if $k=l$).
\end{definition}

\begin{customthm}{B}{\cite[Theorem 4.3]{kirousis2019aggregation}}\label{thm:tot-block}
Let $X$ be a set of feasible evaluations. The following two statements are equivalent.

\begin{enumerate}
\item $X$ is totally blocked.

 \item $X$ admits no binary non-dictatorial aggregator.
\end{enumerate}
\end{customthm}

In the Boolean framework, Theorem \ref{thm:tot-block} can be obtained from Theorem \ref{thm:boolean} above and from \cite[Claim $3.6$]{dokow2010aggregation}. So by applying once more Theorem \ref{thm:boolean}, we get the following result.

\begin{corollary}{Implicit in \cite{dokow2010aggregation}}\label{cor:boolean}
Let $X\subseteq \{0,1\}^m$ be a set of feasible evaluations. The following two statements are equivalent.
\begin{enumerate}
\item $X$ is a possibility domain.
\item $X$ is affine or $X$ admits a binary non-dictatorial aggregator.
\end{enumerate}
\end{corollary}

In the non-Boolean case, Kirousis et al. showed that:
\begin{customthm}{C}{\cite[Theorem 3.1]{kirousis2019aggregation}}\label{thm:pd} Let $X$ be a set of feasible evaluations. The following two statements are equivalent.
\begin{enumerate}
\item  $X$ is a possibility domain.
 \item $X$ admits a binary non-dictatorial aggregator, or a majority aggregator, or a minority aggregator.
     \end{enumerate}
\end{customthm}

We illustrate the two preceding theorems with several examples.
\begin{example}\label{ex:impos}
Let $X_3=\{(1,0,0), (0,1,0), (0,0,1)\}\subseteq\{0,1\}^3$ be the set of all Boolean triples that contain exactly one $1$.  By Theorems \ref{thm:boolean} and \ref{thm:tot-block}, the set $X_3$ is an impossibility domain, since it is not affine ($\oplus((1,0,0), (0,1,0), (0,0,1))=(0,0,0)\notin X_3$) and it does not have a binary non-dictatorial aggregator. For the latter, one has to check each of the $62$ possible $3$-tuples of conservative binary functions over $\{0,1\}$, which is a fairly tedious but straightforward task.
\end{example}

Before proceeding with the next example, let us state a straightforward result that we will need in what follows.
\begin{lemma}\label{lem:cartprod}
Let $Y\subseteq A^l$, $Z\subseteq A^{m-l}$ and assume $G=(g_1,\ldots,g_l)$ and $H=(h_1,\ldots,h_{m-l})$ are $n$-ary aggregators of $Y$ and $Z$ respectively. Define $X\subseteq A^m$ as the \emph{Cartesian product} of $Y$ and $Z$:$$X:=Y\times Z:=\{(x_1,\ldots,x_l,x_{l+1},\ldots,x_m)\mid (x_1,\ldots,x_l)\in Y\text{ and }(x_{l+1},\ldots,x_m)\in Z\}.$$ Then, the $m$-tuple $F=(f_1,\ldots,f_m)$ of $n$-ary operations:\begin{equation*}f_j=\begin{cases}g_j,&\text{if }j=1,\ldots,l,\\ h_{j-l},&\text{if } j=l+1,\ldots,m,
\end{cases}\end{equation*} is an aggregator for $X$.
\end{lemma}
\textbf{Proof}
By the definition of $X$, for any $\ba^i=(a^i_1,\ldots,a^i_m)\in X$, it holds that $(a_1^i,\ldots,a_l^i)\in Y$ and $(a_{l+1}^i,\ldots,a_m^i)\in Z$, $i=1,\ldots,n$. Now, since $G$ and $H$ are aggregators for $Y$ and $Z$, it holds that $G(a_1^i,\ldots,a_l^i):=(b_1,\ldots,b_l)\in Y$ and $H(a_{l+1}^i,\ldots,a_m^i):=(b_{l+1},\ldots,b_m)\in Z$. Thus, again by the definition of $X$, $F(\ba^1,\ldots,\ba^n)=(b_1,\ldots,b_m)\in X$.\hfill$\Box$
\begin{example}
The following two sets $X_4$ and $X_5$  are possibility domains. For $X_4$, we show that this happens using three agents (thus the non-dictatorial aggregator attesting to that is ternary) and, for $X_5$, with two (and thus we find a binary non-dictatorial aggregator).
\begin{itemize}
    \item[(i)] Let $X_4=\{(0,1,2), (1,2,0), (2,0,1), (0,0,0)\}$. This set has been studied in \cite[Example 4]{dokow2010aggregation}. Let $F=(f_1,f_2,f_3)$ be such that, for each $j=1,2,3$:
    \begin{equation*}
	f_j(x, y, z) = \begin{cases}
 {\rm maj}(x,y,z) & \text{ if }  |\{x,y,z\}|\leq 2,\\
 0 & \text{ else}.
 \end{cases}
\end{equation*}
Clearly, $F$ is a majority operation.
To see that $F$ is an aggregator for $X_4$,  we need to check that $F(a,b,c)\in X_4$, only when $a,b,c$ are \emph{pairwise distinct} vectors of $X_4$. In this case, $F(a,b,c)=(0,0,0)\in X_4$, since the input of each $f_j$ contains either two zeros or three pairwise distinct elements.
\item[(ii)] Let  $X_5=X_3\times X_3$, where $X_3$ is as in Example \ref{ex:impos}. It is straightforward to check that $(pr_1^2,pr_1^2,pr_1^2,pr_2^2,pr_2^2,pr_2^2)$ is a non-dictatorial aggregator for $X_5$. In fact, a stronger fact holds: if $Y$ and $Z$ are arbitrary sets, then their Cartesian product $Y\times Z$ is a possibility domain, since it admits non-dictatorial aggregators of any arity $n\geq 2$, defined as the $d$-th projection ${\rm pr}_d^n$ on coordinates from $Y$ and as the $d'$-th projection ${\rm pr}_{d'}^n$ on coordinates from $Z$, where $1\leq d,d'\leq n$ and $d\neq d'$.
\end{itemize}
\end{example}

\subsection{Uniform Possibility Domains}\label{sub:upd}
Let $F=(f_1,\ldots,f_m)$ be an $n$-ary aggregator for $X$.  We say that $F$ is a \emph{uniform non-dictatorial} aggregator for $X$ (of arity $n$) if, for all $j\in \{1,\ldots,m\}$ and for every two-element subset $B_j\subseteq X_j$, it holds that: $$f_j{\restriction_{B_j}}\neq pr_d^n,\text{ for all }d\in \{1,\ldots,n\}.$$ We say that a set $X$ is a \emph{uniform possibility domain} if it admits a uniform non-dictatorial aggregator.

In the non-Boolean framework, uniform non-dictatorial aggregators were introduced by Kirousis et al. \cite{kirousis2019aggregation}. In the Boolean framework, this notion corresponds to \emph{locally non-dictatorial aggregators}, which were introduced by Nehring and Puppe \cite{nehring2010abstract}. Following their terminology, we say that a uniform possibility domain in the Boolean framework is a \emph{local possibility domain}.

The aforementioned sets $X_1$, $X_2$ and $X_4$ are uniform possibility domains, as $X_1$ and $X_4$ admit a majority aggregator, while $X_2$ admits a minority aggregator.
Clearly, if $X$ is a uniform possibility domain, then $X$ is also a possibility domain. The converse, however, is not true.  Indeed, suppose that $X$ is a Cartesian product $X= Y \times Z$, where $Y\subseteq A^l$ and $Z \subseteq A^{m-l}$, with $1 \leq l <m$. If $Y$ or $Z$ is an impossibility domain, then $X$ is not a uniform possibility domain, although it is a possibility domain, since, as seen earlier,  every Cartesian product of two sets is a possibility domain.
 It is also clear that if $Y$ and $Z$ are uniform possibility domains, then so is their Cartesian product $Y\times Z$. In particular, the Cartesian product $X_1\times X_2$ is a uniform possibility domain.

The next result characterizes {\un} possibility domains. It is the stepping stone towards showing  that the following decision problem is solvable in polynomial time:
given a set  $X$ of feasible evaluations, is $X$ a {\un} possibility domain? We first need to define some operators. We say that $f:A^n\to A$  is a \emph{weak near-unanimity  operation (WNU)} (see \cite{larose2017algebra}) if, for all $x,y\in A$, we have that
\begin{equation*}
f(y,x,x,\ldots,x)=f(x,y,x,\ldots,x)=\ldots=f(x,x,x,\ldots,y).
\end{equation*}
In particular, a ternary weak near-unanimity operation is  a function $f:A^3\to A$ such that for all $x,y \in A$, we have that $f(y,x,x) = f(x,y,x)= f(x,x,y).$ Thus, the notion of a ternary weak near-unanimity operation is a common generalization of the notions of
a majority operation and a minority operation.

As with the majority/minority aggregators, we say that $X$ \emph{admits a ternary weak near-unanimity aggregator} $F=(f_1,\ldots,f_m)$, if it admits a ternary aggregator every component of which is a weak near-unanimity operation, i.e., for all $j=1, \ldots, m$
and for  all $x,y \in X_j$, we have that $f_j(y,x,x) = f_j(x,y,x)= f_j(x,x,y).$

Finally, consider the Boolean operations $\wedge,\vee:\{0,1\}^2\mapsto\{0,1\}$, defined as:
\begin{equation*}
	\wedge (x, y) = \begin{cases}
 0 & \text{ if }  x=0\text{ or }y=0,\\
 1 & \text{ else }
 \end{cases}
\end{equation*}
and
\begin{equation*}
	\vee (x, y) = \begin{cases}
 1 & \text{ if }  x=1\text{ or }y=1,\\
 0 & \text{ else. }
 \end{cases}
\end{equation*}
For any $x,y,z\in\{0,1\}$, we define the ternary operations $\wedge^{(3)}$ and $\vee^{(3)}$ as follows:$$\wedge^{(3)}(x,y,z)=\wedge(\wedge(x,y),z)\text{ and }\vee^{(3)}(x,y,z)=\vee(\vee(x,y),z).$$
From now on, to make notation easier, we arbitrarily identify any binary subset $B\subseteq A$ with $\{0,1\}$, in order for the two symmetric operators $\wedge,\vee$ to have meaning on its elements, without needing to assume each time a different ordering for each such $B$.

 \begin{customthm}{D}{\cite[Theorem 5.5]{kirousis2019aggregation}}\label{thm:updcar}
Let $X$ be a set of feasible evaluations. The following three statements are equivalent.

\begin{enumerate}
\item $X$ is a {\un}  possibility domain.
\item $X$ admits a ternary aggregator $F=(f_1,\ldots,f_m)$ such that, for all $j\in\{1,\ldots,m\}$ and for all subsets $B_j\subseteq X_j$ of size $2$, it holds that $f_j{\restriction_{B_j}}\in\{\wedge^{(3)},\vee^{(3)},\maj,\oplus\}.$
 \item $X$ admits a ternary weak near-unanimity aggregator.
\end{enumerate}
\end{customthm}

\subsection{Computational Complexity}\label{sub:PH}
We end this section by briefly discussing the \emph{polynomial hierarchy} (\PH), which
%consists of the most known complexity classes inside \PSPACE, that is, the classes of problems solvable by algorithms needing polynomial space.
consists of a family of complexity classes that contain \NP~and \coNP, and are, in turn, contained in \PSPACE, the class of all decision problems solvable by algorithms that use a polynomial amount of space.

The complexity classes in the polynomial hierarchy were originally defined via \emph{oracles}. An algorithm with access to an oracle $O$ is an algorithm that can, in one step, obtain a yes/no answer for any instance of the decision problem $O$.
For example, an algorithm with access to a \SAT \ oracle, can, in one step, learn if a propositional formula that has come up in its execution, is satisfiable or not. Note that since \SAT \ is \NP-complete, this gives the additional power to the algorithm of deciding any problem in \NP \ in polynomial time.
Given two complexity classes $A$ and $B$, the class of problems that can be decided by an algorithm in $A$, with oracle access to the class $B$ (i.e., with access to some complete problem of $B$), is denoted by $A^B$. In this way, the complexity classes of the \PH~are $\rS_k^P$, $\rP_k^P$ and $\rD_k^P$, $k\in\N$, which are recursively defined as follows:\begin{itemize}
    \item $\rS_0^P=\rS_0^P=\rD_0^P=\P$ and
    \item $\rS_{k+1}^P$ is $\NP$ with oracle $\rS_k^P$, $\rP_{k+1}^P$ is $\coNP$ with oracle $\rS_k^P$ and $\rD_{k+1}^P$ is $\P$ with oracle $\rS_k^P$, $k\in\N$.
\end{itemize}
It is known that
$$\rS_k^P\cup\rP_k^P\subseteq \Delta_{k+1}^P\subseteq \rS_{k+1}^P\cap\rP_{k+1}^P,\text{ }\forall k\in\N.$$ Furthermore, if for some $k\in\N$, we have $\rS_k^P=\rP_k^P$, then \PH \ collapses to that level, in the sense that $\rP_l^P=\rP_l^P$, for all $l\geq k$. For example, if \NP=\coNP, then $\rS_k^P=\rP_k^P=\NP$, for all $k\geq 1$.
Every level of \PH~contains complete problems that are generalizations of SAT. For example, $\rS_k$-\SAT~is complete for $\rS_k^P$, where
$\rS^k$-\SAT~is the following decision problem: given an expression of the form $\exists {\bf x}_1 \forall {\bf x}_2 \ldots Q{\bf x}_k \varphi({\bf x}_1,{\bf x}_2, \ldots, {\bf x}_k)$, where  $\varphi$ is a Boolean formula, is this expression true when the quantifiers vary over the set $\{0,1\}$? (Here, $Q=\exists$ if $k$ is odd, while $Q=\forall$ if $k$ is even.)

 For a more in depth discussion of the polynomial hierarchy, we refer the interested reader to Stockmeyer \cite{stockmeyer1976polynomial}.
Finally, recall that \NP~and~\coNP~can be defined via certificates:\begin{align*}
    \NP= &\{x\in\{0,1\}^*\mid \exists w\in\{0,1\}^{p(|x|)}:\langle x,w\rangle\in \P\},\\
  \coNP= & \{x\in\{0,1\}^*\mid \forall w\in\{0,1\}^{p(|x|)}:\langle x,w\rangle\in \P\},
\end{align*}

Under this notation, $w$ is the certificate, i.e., a candidate solution for the instance $x$, whose validity can be checked in polynomial time. Following this, we obtain an alternative description of the polynomial hierarchy via certificates:
\begin{align*}
    \rS_{k+1}^P= & \exists\rP_k^P=\{x\in\{0,1\}^*\mid \exists w\in\{0,1\}^{p(|x|)}:\langle x,w\rangle\in \rP_k^P\},\\
  \rP_{k+1}^P= & \forall\rS_k^P=\{x\in\{0,1\}^*\mid \forall w\in\{0,1\}^{p(|x|)}:\langle x,w\rangle\in \rS_k^P\},
\end{align*}
where $p(|x|)$ is polynomial in the size of $|x|$. We will use both definitions of these complexity classes to derive our results.

\section{Explicitly given Domains}\label{sec:abstract}

In subsection \ref{ssec:poss-tract}, we present a polynomial-time algorithm for checking if $X$ is a possibility domain. In subsection \ref{ssec:tcs}, we prove that this problem is also expressible in Transitive Closure Logic (TCL). Finally, in subsection \ref{ssec:upd-tract}, we provide a polynomial-time algorithm for checking if a domain $X$ is a uniform possibility domain.

\subsection{Tractability of Possibility Domains} \label{ssec:poss-tract}

Theorems \ref{thm:boolean} and \ref{thm:pd} provide necessary and sufficient conditions for a set $X$ to be a possibility domain in the Boolean framework and in the non-Boolean framework, respectively. Admitting a binary non-dictatorial aggregator is a condition that appears in both of these characterizations. Our first result asserts that this condition can be checked in polynomial time.

\begin{theorem}\label{thm:binary-tract}
There is a polynomial-time algorithm for solving the following problem: given a set $X$ of feasible evaluations, determine whether or not $X$ admits a binary non-dictatorial aggregator and, if it does, produce one.
\end{theorem}

We first show that the existence of a binary non-dictatorial aggregator on $X$ is tightly related to connectivity properties of a certain directed graph $H_X$ defined next. If $X\subseteq A^m$ is a set of feasible evaluations, then $H_X$ is the following directed graph:
\begin{itemize}
\item The  vertices of $H_X$ are the pairs of \emph{distinct} elements $u,u'\in X_j$, for  $j\in \m$.
Each such vertex will usually be denoted by $uu'_j$. When the coordinate $j$ is understood from the context, we will often  be dropping the subscript $j$, thus denoting such a vertex by $uu'$.

Also, if $u \in X_j$, for some $j\in \m$, we will often use the notation $u_j$ to indicate that $u$ is an element of $X_j$.

\item Two vertices $uu'_k$ and $vv'_l$, where $k\neq l$, are connected by a directed  edge from $uu'_k$  to $vv'_l$,  denoted by $uu'_k \rightarrow vv'_l$, if there are a total evaluation $z \in X $ that extends the partial evaluation $(u_k,v_l)$ and a total evaluation $z' \in X $ that extends the partial evaluation  $(u'_k,v'_l)$, such   that there is no  total evaluation $y \in X $  that extends $(u_k,v'_l)$, and has the property that $y_i = z_i$ or $y_i =z'_i$,  for every $i \in \{1, \ldots, m\}$.
\end{itemize}
For vertices $uu'_k$, $vv'_l$, corresponding to  issues $k, l$ (that need not be distinct),  we
 write $uu'_k \rightarrow \rightarrow vv'_l$
to  denote the existence of a directed path from $uu'_k$   to $vv'_l$. Note that if $uu'_k\rightarrow vv'_l$, then $v'v_l\rightarrow u'u_k$, since we require for the same partial vectors to be extendable.

In the next example, we describe explicitly the graph $H_X$ for several different sets $X$ of feasible voting patters. Recall that a \emph{directed} graph $G$ is  \emph{strongly connected} if for every pair of vertices $(u,v)$ of $G$, there is a (directed) path from $u$ to $v$.

Finally, consider the graph $G_X$ of Definition \ref{def:tot-blockr}. As will become apparent in the sequel, $H_X$ is strongly connected if and only if $G_X$ is. Nevertheless, $G_X$ can have more edges, since the infeasible partial evaluation $x$ can be changed in any $j\in K$ in order to become feasible.

\begin{example}\label{ex:HX}
Recall the two Boolean domains $X_2=\{(0,0,1),(0,1,0),(1,0,0),(1,1,1)\}$ and $X_3=\{(0,0,1),(0,1,0),(1,0,0)\}$ of Examples  \ref{ex:bijaff} and \ref{ex:impos}. Both $H_{X_2}$ and $H_{X_3}$ have six vertices, namely $01_j$ and $10_j$, for $j=1,2,3$. In the figures below, we use undirected edges between two vertices $uu'_k$ and $vv'_l$ to denote the existence of both $uu'_k \rightarrow vv'_l$ and $vv'_l \rightarrow uu'_k$.

\begin{figure}[H]
\begin{center}
\begin{tikzpicture}[node/.style={circle,draw=black!100,fill=white!20,minimum size=18pt,inner sep=0pt},nonode/.style={circle}]

\draw (0,-2) -- (0,2);

%supporting nodes
\node[nonode] (c) {};
\node[nonode] (g) [above = 2cm of c] {};

\node[nonode] (u) [above = of c] {};
\node[nonode] (d) [below = of c] {};

\node[nonode] (cu) [above = 0.2cm of c] {};
\node[nonode] (cd) [below = 0.2cm of c] {};

\node[nonode] (bu) [above = 0.3cm of c] {};

%X2: u1 --> 01_1, u2 --> 10_1, u3 --> 01_2, u4 --> 10_2
%u5--> 01_3, u_6 --> 10_3

\node[nonode] (X2) [left=2.375cm of g] {$H_{X_2}$};

\node[node] (u1) [left=5cm of u] {$01_1$};
\node[node] (u2) [left=5cm of bu] {$10_1$};
\node[node] (u3) [left=0.25cm of u] {$01_2$};
\node[node] (u4) [left=0.25cm of bu] {$10_2$};
\node[node] (u5) [left=3.125cm of d] {$01_3$};
\node[node] (u6) [left=2.125cm of d] {$10_3$};

\draw (u1) -- (u3);
\draw (u1) -- (u4);
\draw (u1) -- (u5);
\draw (u1) -- (u6);
\draw (u2) -- (u3);
\draw (u2) -- (u4);
\draw (u2) -- (u5);
\draw (u2) -- (u6);
\draw (u3) -- (u5);
\draw (u3) -- (u6);
\draw (u4) -- (u5);
\draw (u4) -- (u6);

%X3: v1 --> 01_1, v2 --> 10_2, v3 --> 01_3, v4 --> 10_1
%v5--> 01_3, v_6 --> 10_3

\node[nonode] (X3) [right=2.375cm of g] {$H_{X_3}$};

\node[node] (v1) [right=0.25cm of cu] {$01_1$};
\node[node] (v2) [right=2.6cm of u] {$10_2$};
\node[node] (v3) [right=5cm of cu] {$01_3$};
\node[node] (v4) [right=5cm of cd] {$10_1$};
\node[node] (v5) [right=2.6cm of d] {$01_2$};
\node[node] (v6) [right=0.25cm of cd] {$10_3$};

\draw (v1) -- (v2);
\draw (v2) -- (v3);
\draw (v3) -- (v4);
\draw (v4) -- (v5);
\draw (v5) -- (v6);
\draw (v6) -- (v1);

\end{tikzpicture}
\end{center}
\end{figure}

Consider $01_1$, $01_2$ of $H_{X_2}$. Since the partial vectors on the first two coordinates, $(0,0)$ and $(1,1)$, extend to $(0,0,1)$ and $(1,1,1)$, respectively, we need to check if there is a vector in $X_2$ extending $(0,1)$, but whose third coordinate is $1$. Since $(0,1,1)\notin X_2$, we have that $H_{X_2}$ contains both edges $01_1\rightarrow 01_2$ and $01_2\rightarrow 01_1$. Now, since the partial vectors, again in the first two coordinates, $(0,1)$ and $(1,0)$ extend to $(0,1,0)$ and $(1,0,0)$, respectively, and since neither $(0,0,0)$ nor $(1,1,0)$ are in $X_2$, we have that $01_1\leftrightarrow 10_2$. By the above and because of the symmetric structure of $X_2$, it is easy to see that every two vertices $uu'_i$ and $vv'_j$ of $H_{X_2}$ are connected if and only if $i\neq j$.

For $X_3$, observe that, since \emph{no} partial vector containing two ``$1$"'s, in any two positions,  extends to an element of $X_3$, there are no edges between the vertices $01_i$, $01_j$ and $10_i$, $10_j$, for any $i,j\in\{1,2,3\}$, $i\neq j$. In the same way as with $H_{X_2}$, we get that $H_{X_3}$ is a cycle.

There are two observations to be made, concerning $H_{X_2}$ and $H_{X_3}$. First, they are both strongly connected graphs. Also, neither $X_2$ nor $X_3$ admit binary non-dictatorial aggregators. $X_2$ admits only a minority aggregator (the proof is left for the interested reader) and $X_3$ is an impossibility domain (as shown in Example \ref{ex:impos}).

Finally, consider $X_6:=\{(0,1),(1,0)\}$. The graph $H_{X_6}$ has four vertices, $01_1$, $10_1$, $01_2$ and $10_2$, and it is easy to see that $H_{X_6}$ has only the following edges:

\begin{figure}
\begin{center}
\begin{tikzpicture}[node/.style={circle,draw=black!100,fill=white!20,minimum size=18pt,inner sep=0pt},nonode/.style={circle}]

\node[nonode] (X6) {$H_{X_6}$};
\node[nonode] (h) [below =0.25cm of X6] {};

\node[node] (w1) [left =0.5cm of h] {$01_1$};
\node[node] (w2) [right =0.5cm of h] {$10_2$};
\node[node] (w3) [below =0.25cm of w1] {$10_1$};
\node[node] (w4) [below =0.25cm of w2] {$01_2$};

\draw (w1) -- (w2);
\draw (w3) -- (w4);

\end{tikzpicture}
\end{center}
\end{figure}

Observe that $X_6$ is not strongly connected (it is not even connected) and that,  in contrast to the sets $X_2$ and $X_3$, the set $X_6$ admits two binary non-dictatorial aggregators, namely, $(\wedge,\vee)$ and $(\vee,\wedge)$. In Lemma \ref{lem:binary}, we establish a tight connection between strong connectedness and the existence of binary non-dictatorial aggregators.

\end{example}

We now state and prove two lemmas about the graph $H_X$.

\begin{lemma}\label{lem:propag}
Assume that $F = (f_1, \ldots, f_m)$ is  a binary aggregator on $X$.
\begin{enumerate}
\item If   $uu'_k \rightarrow vv'_l$ and  $f_k(u,u')= u$, then $f_l(v,v') = v$.
% moreover,
%if  $f_l(v,v')= v'$, then $f_k(u,u') = u'$.
\item
If  $uu'_k \rightarrow\rightarrow  vv'_l$ and $f_k(u,u')= u$, then $f_l(v,v') = v$.
% moreover,
%if  $f_l(v,v')= v'$, then $f_k(u,u') = u'$.
\end{enumerate}
\end{lemma}
\textbf{Proof}
The first part of the lemma follows from the definitions and the fact that  $F$ is conservative. Indeed, if $uu'_k \rightarrow vv'_l$, then there are a total evaluation $z=(z_1,\ldots,z_m) \in X$ that extends $(u_k,v_l)$ (i.e., $z_k=u$ and $z_l=v$)  and a total evaluation $z'=(z'_1,\ldots,z'_m)\in X$  that extends $(u'_k,v'_l)$ (i.e., $z'_k=u'$ and $z'_l=v'$), such that there is no total evaluation in $X$ that extends $(u_k,v'_l)$ and agrees with $z$ or with $z'$ on every coordinate. Consider the total evaluation $(f_1(z_1,z'_1), \dots, f_m(z_m,z'_m))$, which is in $X$ because $F$ is an aggregator on $X$. Since each $f_j$ is conservative, we must have that $f_j(z_j,z'_j) \in \{z_j,z'_j\}$, for every $j$, hence $f_l(z_l,z'_l)=f_l(v,v') \in \{v,v'\}$. Consequently, if $f_k(u,u')=u$, then we must have $f_l(v,v') = v$, else $(f_1(z_1,z'_1), \dots, f_m(z_m,z'_m))$ extends $(u_k,v'_l)$ and agrees with $z$ or with $z'$ on every coordinate. The second part of the lemma easily follows from the first part by induction.\hfill$\Box$

\begin{lemma}\label{lem:binary}
$X$ admits a binary non-dictatorial aggregator if and only if the directed graph $H_X$   is not strongly connected.
\end{lemma}
Before delving into the proof, consider the graphs of Example \ref{ex:HX}. Using the fact that the graphs $H_{X_2}$ and $H_{X_3}$ are strongly connected and also using the second item of Lemma \ref{lem:propag}, it is easy to see that $X_2$ and $X_3$ admit no binary non-dictatorial aggregator; indeed,  let $F=(f_1,f_2,f_3)$ be a binary aggregator of either of  these two sets and suppose that $f_1(0,1)=0$. Since in both graphs $H_{X_2}$ and $H_{X_3}$, there are paths from $01_1$ to \emph{every} other vertex, it follows that $f_j=pr_1^2$, $j=1,2,3$. If $f_1(0,1)=1$, we get that $f_j=pr_2^2$, $j=1,2,3$, in the same way.

In contrast, consider $H_{X_6}$ and
  let $G=(g_1,g_2)$  be a pair of binary functions with $g_1(0,1)=0$. For $G$ to be an aggregator, Lemma \ref{lem:propag} forces us to set $g_2(1,0)=1$. But now, by setting $g_1(1,0)=0$, and thus $g_2(0,1)=1$, we get that $(g_1,g_2)=(\wedge,\vee)$ is a non-dictatorial binary aggregator for $X_6$.\vspace{0.2cm}

  \noindent\textbf{Proof of Lemma \ref{lem:binary}} We first show that if $X$ admits a binary non-dictatorial aggregator, then $H_X$ is not strongly connected. In the contrapositive form, we show that if $H_X$ is strongly connected, then $X$ admits no binary non-dictatorial aggregator. This is an easy consequence of the preceding Lemma \ref{lem:propag}. Indeed, assume that $H_X$ is strongly connected and let $F=(f_1,\ldots,f_m)$ be a binary aggregator on $X$. Take two distinct elements $x$ and $x'$ of $X_1$. Since $F$ is conservative, we have that $f_1(x,x')\in \{x,x'\}$. Assume first that $f_1(x,x')=x$. We claim that $f_j= {\rm pr}^2_1$, for every $j\in \m$. To see this, let $y$ and $y'$ be two distinct elements of $X_j$, for some $j\in \m$. Since $H_X$ is strongly connected, we have that $xx'_1 \rightarrow\rightarrow yy'_j$. Since also $f_1(x,x')=x$, Lemma \ref{lem:propag} implies that $f_j(y,y')=y = {\rm pr}^2_1(y,y')$ and so $f_j = {\rm pr}^2_1$.
 Next, assume that $f_1(x,x') = x'$. We claim that $f_j= {\rm pr}^2_2$, for every $j\in \m$. To see this, let $y$ and $y'$ be two distinct elements of $X_j$, for some $j\in \m$. Since $H_X$ is strongly connected, we have that $yy'_j\rightarrow\rightarrow xx'_1$, hence, if $f_j(y,y')=y$, then, Lemma \ref{lem:propag}, implies that $f_1(x,x') = x$, which is a contradiction because $x\not = x'$. Thus, $f_j(y,y')=y'$ and so $f_j = {\rm pr}^2_2$.

For the converse, assume that $H_X$ is not strongly connected and let $uu'_k$, $vv'_l$ be two vertices of $H_X$ such that there is no path from $uu'_k$ to $vv'_l$ in $H_X$, i.e.,  it is not true that $uu'_k \rightarrow \rightarrow vv'_l$. Let $V_1, V_2$ be a partition of the vertex set such that $uu'_k \in V_1, vv'_l \in V_2$, and there is  no edge from a vertex in $V_1$ to a vertex in $V_2$. We will now define a binary aggregator $F = (f_1, \ldots, f_m)$ and prove that it is non-dictatorial.
%, which will complete the proof.

Given $z,z' \in X $,
we set  $f_j(z_j, z'_j) = z_j$ if $zz'_j \in V_1$, and we set
$f_j(z_j, z'_j) = z'_j$ if $zz'_j \in V_2$, for $j\in \m$. Since $uu'_k \in V_1$, we have that $f_k\not = {\rm pr}^2_2$; similarly, since
$vv'_l\in V_2$, we have that  $f_l\not = {\rm pr}^2_1$. Consequently, $F$ is not a dictatorial function on $X$.
Thus, what remains to be proved is that if $z,z,'\in X$, then  $F (z, z') \in X$. For this, we will show that if $ F (z, z')  \not\in X$, then there is an edge from an element of  $V_1$ to an element of $V_2$, which is a contradiction.

Assume that $q = F (z, z') \not \in X$. Let $K$ be a minimal subset of $\{1, \ldots, m\}$ such that $ q{\restriction K} $ cannot be extended to a total evaluation $w$ in  $X$ that agrees with $z$ or with $z'$ on  $\{1, \ldots, m\}\setminus K$ (i.e., if $j\in \{1, \ldots, m\}\setminus K$, then $w_j=z_j$ or $w_j= z'_j$).
 Since $z'$ is in $X$, it does not extend   $q{\restriction K} $, hence there is a number $s\in K$ such that $q_s=f_s(z_s,z'_s)= z_s \not = z'_s$. It follows that
 the vertex $zz'_s$ is in $V_1$.
 Similarly, since $z$ is in $X$, it does not extend  $q{\restriction K} $, hence there is a number $t\in K$ such that $q_t=f_t(z_t,z'_t)= z'_t \not = z_t$. It follows that the vertex $zz'_t$ is in $V_2$. Consequently, there is no edge from $zz'_s$ to $zz'_t$ in $H_X$.  We will arrive at a contradiction by showing that $zz'_s \rightarrow zz'_t$, i.e., there is an edge $zz'_s$ to $zz'_t$ in $H_X$. Consider the set $K\setminus \{t\}$. By the minimality of $K$, there is a total evaluation $w$ in $X$ that extends $q{\restriction K\setminus \{t\}}$ and agrees with $z$ or with $z'$  outside $K\setminus \{t\}$. In particular, we have that $w_s=q_s=z_s$ and $w_t=z_t$.  Similarly, by considering the set $K\setminus \{s\}$, we find that there is a  total evaluation $w'$ in $X$ that extends $q{\restriction K\setminus \{s\}}$ and agrees with $z$ or with $z'$  outside $K\setminus \{s\}$. In particular, we have that $w'_s=z'_s$ and $w_t=q_t=z'_t$. Note that $w$ and $w'$ agree on $K\setminus \{s,t\}$. Since $q{\restriction K}$ does not extend to a total evaluation that agrees with $z$ or with $z'$ outside $K$, we conclude that there is no total evaluation $y$ in $X$
   that extends $(z_s,z'_t)$ and agrees with $w$ or with $w'$ on every coordinate. Consequently, $zz'_s \rightarrow zz'_t$, thus we have arrived at a contradiction.\hfill $\Box$\vspace{0.2cm}

\noindent\textbf{Proof of Theorem \ref{thm:binary-tract}:} Given a set $X$ of feasible evaluations, the graph $H_X$ can be constructed in time bounded by a polynomial in the size $|X|$ of $X$ (in fact, in time $O(|X|^5)$.
There are well-known polynomial-time algorithms for testing if a graph is strongly connected and, in case it is not, producing the \emph{strongly connected components (\scc)} of the graph; e.g., Kosaraju's algorithm presented in \cite{sharir1981strong} and \cite{Tarjan72}. Consequently, by Lemma \ref{lem:binary}, there is a polynomial-time algorithm for determining whether or not a given set $X$ admits a binary non-dictatorial aggregator. Moreover, if $X$ admits such an aggregator, then one can be constructed in polynomial-time from the strongly connected components of $H_X$  via the construction in the proof of Lemma \ref{lem:binary}.\hfill$\Box$\vspace{0.2cm}

The next corollary follows   from Theorem \ref{thm:binary-tract} and Theorem \ref{thm:tot-block}.

\begin{corollary} \label{cor:tot-block}
There is a polynomial-time algorithm for the following decision problem:
given a set $X$ of feasible evaluations, is $X$ totally blocked?
\end{corollary}

We now turn to the problem of detecting possibility domains in the non-Boolean framework.

\begin{theorem}\label{thm:pdmain}
There is a polynomial-time algorithm for solving the following problem:  given a set  $X$  of feasible evaluations, determine whether or not $X$ is a possibility domain and, if it is,  produce a binary non-dictatorial aggregator, or a ternary majority aggregator or a ternary minority aggregator for $X$.
\end{theorem}
\textbf{Proof} It is straightforward to check that, by Theorem \ref{thm:pd} and Theorem \ref{thm:binary-tract}, it suffices to show that there is a polynomial-time algorithm that, given $X$, detects  whether or not $X$
admits a majority aggregator or a minority aggregator, and, if it does, produces such an aggregator.

Let $X$ be a set of feasible evaluations, where  $ I = \{1, \ldots, m\}$ is the set of issues and $A$ is the set of the position values. We define the {\em disjoint union} $\A$ of the  set of position values as:
\[
\centerline{$\mathcal{A} = \bigsqcup_{j=1}^m A  = \bigcup_{j=1}^m \{(x,j) \mid x \in A\}$.}\]

\noindent We also set
\[\centerline{
 $\tilde{X} = \{((x_1, 1), \ldots, (x_m,m))  \mid (x_1, \ldots, x_m) \in X \} \subseteq \mathcal{A}^m .$}\]

We will show that we can go back-and-forth between conservative majority or minority polymorphisms for $\tilde{X}$ and majority or minority  aggregators for $X$.

Let $f:\mathcal{A}^n\to \mathcal{A}$ be a conservative polymorphism for $\tilde{X}$. We define the $m$-tuple $F=(f_1,\ldots,f_m)$ of $n$-ary functions $f_1,\ldots,f_m$ as follows:
 if $x_j^1,\ldots,x_j^n\in X_j$, for $j\in \m$, then we set  $f_j(x_j^1,\ldots,x_j^n)=y_j,$
  where $y_j$ is such that
   $f((x_j^1,j),\ldots,(x_j^n,j))=(y_j,j)$. Such a  $y_j$ exists and is one of the $x_j^i$'s because $f$ is conservative, and hence $f((x_j^1,j),\ldots,(x_j^n,j)) \in \{(x_j^1,j),\ldots,(x_j^n,j)\}$.
   It is easy to see that $F$ is an aggregator for $X$. Moreover, in case $n=3$, if $f$ is a majority  or a minority operation on $\tilde{X}$, then $F$ is a majority or a minority aggregator on $X$.

Next, let $F=(f_1,\ldots,f_m)$ be a (ternary) majority or a minority aggregator for $X$.  We define a ternary function $f:\mathcal{A}^3\rightarrow \mathcal{A}$ as follows. Let $(x,j), (y,k), (z,l)$ be three elements of $\mathcal{A}$.
\begin{itemize}
\item If $j=k=l$, then we set $f((x,j), (y,k), (z,l)) = (f_j(x,y,z),j)$.
\item If $j, k, l$ are \emph{not} all equal, then if at least two of $(x,j), (y,k), (z,l)$ are equal to each other,
we set

\centerline{$f((x,j), (y,k), (z,l)) = {\rm maj}((x,j), (y,k), (z,l)),$}

\noindent if $F$ is a majority aggregator on $X$, and we set

\centerline{$f((x,j), (y,k), (z,l)) = \oplus ((x,j), (y,k), (z,l)),$}

\noindent if $F$ is a minority aggregator on $X$;
\item otherwise, we set $f((x,j), (y,k), (z,l)) = (x,j)$.
\end{itemize}
It is easy to see that if $F$ is a majority or a minority aggregator for $X$, then $f$ is a conservative majority or a conservative minority polymorphism on $\tilde{X}$.
It follows that $X$ admits a majority or a minority aggregator if and only if $\tilde{X}$ is closed under a conservative majority or minority polymorphism.
 Bessiere et al. \cite[Theorem 1]{bessiere2013detecting} and Carbonnel \cite[Theorem 1]{carbonnel2016meta} design polynomial-time algorithms that detect if a given constraint language $\Gamma$  has a conservative majority or a conservative minority polymorphism, respectively, and,  when it has, compute such a polymorphism.
Here, we  apply these results to $\Gamma=\{\tilde{X}\}$.\hfill$\Box$\vspace{0.2cm}

We end this subsection by showing that using the graph $H_X$, we can compute a binary aggregator for $X$ that has as many symmetric components as possible. This will allow us to obtain better complexity bounds in Section \ref{sec:implicit}.

Given a domain $X\subseteq A^m$ and a binary aggregator $F=(f_1,\ldots,f_m)$ for $X$, we say that $F$ is a \emph{maximum symmetric} aggregator for $X$ if, for every binary aggregator $G=(g_1,\ldots,g_m)$ for $X$, for every $j\in\m$ and for all binary $B_j\subseteq X_j$, if $g_j{\restriction_{B_j}}$ is symmetric, then so is $f_j{\restriction_{B_j}}$. Note that a maximum symmetric aggregator does not necessarily have any symmetric components, for example in case $X$ is an impossibility domain. Furthermore, if $F$ and $G$ are both maximum symmetric aggregators for $X$, then they necessarily have the same symmetric components.
\begin{lemma}\label{lem:maxsym}
Every domain $X$ admits a maximum symmetric aggregator.\end{lemma}
\textbf{Proof} Assume that there is no maximum symmetric aggregator for $X$. Then, there exist two indices $i,j\in\m$ and two binary subsets $B_i\subseteq X_i$, $B_j\subseteq X_j$, such that:\begin{itemize}
    \item there are two binary aggregators $F=(f_1,\ldots,f_m)$ and $G=(g_1,\ldots,g_m)$ such that $f_i{\restriction_{B_i}}$ and $g_j{\restriction_{B_j}}$ are symmetric and
    \item there is no binary aggregator $H=(h_1,\ldots,h_m)$ such that both $h_i{\restriction_{B_i}}$ and $h_j{\restriction_{B_j}}$ are symmetric.
\end{itemize} Let $H=(h_1,\ldots,h_m)$ be the $m$-tuple of binary functions such that, for all $k\in\m$ and for all $x,y\in X_k$: $h_k(x,y):=g_k(f_k(x,y),f_k(y,x))$. That $H$ is an aggregator for $X$ follows easily from the fact that $F$ and $G$ are. Let also $B_i=\{a,b\}$ and $B_j=\{c,d\}$. Since $f_i{\restriction_{B_i}}$ is symmetric, $f_i(a,b)=f_i(b,a)$, thus:$$h_i(a,b)=g_i(f_i(a,b),f_i(b,a))=g_i(f_i(a,b),f_i(a,b))=f_i(a,b).$$ It follows that $h_i{\restriction_{B_i}}$ is symmetric. Now, by the hypothesis, $f_j{\restriction_{B_j}}$ is not symmetric. Assume w.l.o.g. that $f_j{\restriction_{B_j}}=\pr_1^2$. Thus, it holds that:$$h_j(c,d)=g_i(f_i(c,d),f_i(d,c))=g_i(c,d), $$ which means that $h_j{\restriction_{B_j}}$ is symmetric. Contradiction. \hfill$\Box$\vspace{0.2cm}

\begin{remark}
Note that in Lemma \ref{lem:maxsym} we in fact show that given two binary aggregators $F$ and $G$ that are symmetric on some binary $B_i\subseteq X_i$ and $B_j\subseteq X_j$, then we can always construct a binary aggregator $H$ that is symmetric in both $B_i$ and $B_j$.
\end{remark}

To proceed, we discuss a result concerning the structure of the graph $H_X$. We say that two \scc's $S_p$ and $S_q$ of $H_X$ are \emph{related} if there exists a $j\in\m$ and two distinct elements $u,u'\in X_j$, such that $uu'_j\in S_p$ and $u'u_j\in S_q$.
\begin{lemma}\label{lem:sccHX}
Let $X$ be a set of feasible voting patterns and assume that $S_p$, $S_q$ and $S_r$ are three pairwise distinct \scc's of $H_X$. Then, $S_p$, $S_q$ and $S_r$ cannot be pairwise related.
\end{lemma}
\textbf{Proof}
To obtain a contradiction, assume they are. Then, there exist (not necessarily distinct) indices $i,j,k\in\m$ and pairwise distinct elements $u,u'\in X_i$, $v,v'\in X_j$ and $w,w'\in X_k$ such that $uu'_i$, $vv'_j\in S_p$, $ww'_k$, $u'u_i\in S_q$ and $v'v_j$, $w'w_k\in S_r$. Since $uu'_i\rightarrow\rightarrow vv'_j$ and $vv'_j\rightarrow\rightarrow uu'_i$,  it follows that $v'v_j\rightarrow\rightarrow u'u_i$ and $u'u_i\rightarrow\rightarrow v'v_j$. Thus, $S_q$ and $S_r$ form together an \scc \ of $H_X$. Contradiction.\hfill$\Box$\vspace{0.2cm}

We are now ready to show that we can find a maximum symmetric aggregator for a domain $X$, in polynomial time in its size.
\begin{corollary}\label{cor:maxsym}
There is a polynomial-time algorithm for solving the following problem:  given a set  $X$  of feasible evaluations, produce a maximum symmetric aggregator for $X$.
\end{corollary}
\textbf{Proof} Construct the graph $H_X$. For a set of vertices $S$, let $N^+(S)$ and $N^-(S)$ be its \emph{extended outwards and inwards neighborhood} respectively in $H_X$. That is, $N^+(S)=S\cup\{uu'_i\mid \exists vv'_j\in S:\ vv'_j\rightarrow\rightarrow uu'_i\}$ and $N^-(S)=S\cup\{uu'_i\mid \exists vv'_j\in S:\ uu'_i\rightarrow\rightarrow vv'_i\}$.

We define the $m$-tuple of binary functions $F=(f_1,\ldots,f_m)$ as follows. If $H_X$ is strongly connected, set $f_j=\pr_1^2$ for all $j\in\m$. Else, assume w.l.o.g. that $H_X$ is \emph{connected}. If it is not, we can deal with each connected component independently in the same way. Assume that $S_1,\ldots,S_t$, $t\geq 2$, are the \scc's of $H_X$, in their topological order. \begin{enumerate}
    \item For each $uu'_i\in S_1$, set $f_i(u,u')=u'$.
    \item Let $S$ be the set of vertices of every \scc \ of $H_X$ that is related with $S_1$. For each $vv'_j\in N^+(S)$, set $f_j(v,v')=v$.
    \item Let $S'$ be the set of vertices of every \scc \ of $H_X$ that is related with an $\scc$ of $S$. Note that due to Lemma \ref{lem:sccHX}, such an \scc \ cannot be related with $S_1$. For each $ww'_j\in N^-(S)$, set $f_j(w,w')=w'$.
\end{enumerate}
Note that any remaining \scc \ must be in another connected component. If this is the case, we proceed as above for each such connected component. Finally, to be formally correct, let $f_j(a,b)=a$ for every $j\in\m$ and $a,b\in A$ such that either $a$ or $b\notin X_j$.

We have already shown that, given $X$, $H_X$ can be constructed in polynomial time in its size and its \scc's can also be computed in linear time to the size of $H_X$. Steps $1--3$ can easily be implemented by checking once every \scc \ of $H_X$. Thus, the overall process is clearly polynomial.

It remains to show that $F=(f_1,\ldots,f_m)$ is indeed a maximum symmetric aggregator for $X$. To obtain a contradiction, suppose it is not. Then, there exist a $i\in\m$, a binary subset $B_i\subseteq X_i$ and a binary aggregator $G=(g_1,\ldots,g_m)$ for $X$ such that $f_i{\restriction_{B_i}}$ is not symmetric, whereas $g_i{\restriction_{B_i}}$ is.

Assume that $B_i=\{u,u'\}$. Since $G$ is a binary aggregator for $X$ and $g_i{\restriction_{B_i}}$ is symmetric, by Lemma \ref{lem:propag} there are no paths $uu'_i\rightarrow\rightarrow u'u_i$ or $u'u_i\rightarrow\rightarrow uu'_i$ in $H_X$. Consequently, there are two distinct \scc's of $H_X$, say $S_p$ and $S_q$, such that $uu'_i\in S_p$ and $u'u_i\in S_q$ and there are no paths connecting a vertex in $S_p$ with a vertex in $S_q$. Given the way we constructed $F=(f_1,\ldots,f_m)$, the vertices of both these \scc's are either all in $S$ or they are all in $S'$. We show that in both cases, there exist three pairwise distinct \scc's of $H_X$ that are pairwise related. This is a contradiction, by Lemma \ref{lem:sccHX}.

First, assume that the vertices of $S_p$ and $S_q$ are all in $S$. The case where the vertices of $S_p$ and $S_q$ are all in $S'$ is analogous. If both of them are related with $S_1$, then $S_1$, $S_p$ and $S_q$ are pairwise related. Contradiction. Else, without loss of generality, assume that $S_p$ is not related to $S_1$. Then, there exists some \scc \ $S_r$ of $H_X$ that is related with $S_1$, such that $S_p\subseteq N^+(S_r)$. Since $S_r$ is related with $S_1$, there is some vertex $vv'_j\in S_1$ such that $v'v_j\in S_r$. Then $v'v_j\rightarrow\rightarrow uu'_i$, which implies that $u'u_i\rightarrow\rightarrow vv'_j$. Contradiction, since we took the \scc's of $H_X$ in their topological order. \hfill$\Box$

A direct implication of Corollary \ref{cor:maxsym}, that will be used later, is the following.
\begin{corollary}\label{HXmaxsym}
%Given the graph $H_X$, there is a polynomial-time procedure that produces an maximum symmetric aggregator $F$ for $X$.
There is a polynomial-time procedure that, given the graph $H_X$,  the procedure produces a maximum symmetric aggregator $F$ for $X$.
\end{corollary}

\subsection{Expressibility in Transitive Closure Logic}\label{ssec:tcs}
We now establish that Theorem \ref{thm:binary-tract} can be refined to show that testing whether $X$ admits a binary non-dictatorial aggregator can be expressed in Transitive Closure Logic. Before spelling out the technical details, we present a minimum amount of the necessary background from mathematical logic and descriptive complexity.

\paragraph{First-Order Logic, Least Fixed-Point Logic, and Transitive Closure Logic}
First-order logic is a formalism for specifying properties of mathematical objects, such as graphs, trees, partial orders, and, more generally, relational structures. By definition, a \emph{relational schema} is a tuple ${\bf R}=(R'_1,\ldots,R'_m)$ of relational symbols $R'_i$, $1\leq i\leq m$, each of which has a specified natural number $r_i$ as its arity. A \emph{relational structure} over such a schema $\bf R$ is a tuple of the form ${\bf A}=(A,R_1,\ldots, R_m)$, where $A$ is a set called the \emph{universe} of $\bf A$ and each $R_i$ is a relation on $A$ of arity $r_i$, $1\leq i\leq m$. For example, a graph is a relational structure of the form ${\bf G}=(V,E)$, where $E$ is a binary relation on $V$. A \emph{finite} relational structure is a relational structure with a finite set as its universe.
A formula of \emph{first-order logic} is an expression built from \emph{atomic} formulas of the form $x_i=x_j$ and $R_k(x_1,\ldots,x_{r_k})$ using conjunctions, disjunctions, negations, universal, and existential quantification. The semantics of formulas of first-order logic are given by interpreting the quantifiers $\exists$ and $\forall$ as ranging over the universe of the relational structure at hand. For example, in the case of graphs, the first-order formula $\forall x\forall y(E(x,y)\lor \exists z(E(x,z)\land E(z,y))$ asserts that the graph has diameter at most $2$. For the precise definition of the syntax and the semantics of first-order logic, we refer the reader to \cite{enderton2001mathematical}.

It is well known that first-order logic has rather limited expressive power on finite structures. In particular, there is no formula of first-order logic that expresses \emph{connectivity} on finite graphs; this means that there is no formula $\psi$ of first-order logic such that a finite graph $\bf G$ satisfies $\psi$ if and only if $\bf G$ is connected. Moreover, the same holds true for other properties of finite graphs of algorithmic significance, such as \emph{acyclicity} and  $2$-\emph{colorability}; for details, see, e.g., \cite{DBLP:books/sp/Libkin04}.
Intuitively, the reason for these limitations of first-order logic is that first-order logic on finite structures lacks a recursion mechanism.

Least Fixed-Point Logic (LFP) augments first-order logic with a recursion mechanism in the form of least fixed-points of \emph{positive} first-order formulas. More formally, one considers first-order formulas of the form $\varphi(x_1,\ldots,x_n,S)$, where $\varphi(x_1,\ldots,x_n,S)$ is a first-order formula over a relational schema with an extra $n$-ary relation symbol $S$ such that every occurrence of $S$ is within an even number of negation symbols. Every such formula has a \emph{least fixed-point}, that is, for every relational structure $\bf A$, there is a smallest relation $S^*$ such that $S^*=\{(a_1,\ldots,a_n)\in A^n: {\bf A}\models \varphi(a_1,\ldots,a_n,S^*)\}$.
We use the notation $\varphi^\infty(x,y)$ to denote a new formula that expresses the least fixed-point of $\varphi(x_1,\ldots,x_n,S)$.
 For example, if $\varphi(x_1,x_2,S)$ is the formula $E(x_1,x_2) \lor \exists z (E(x_1,z)\land S(z,x_2))$, then, for every graph ${\bf G}=(V,E)$, the least fixed-point $\varphi^\infty(x_1,x_2)$ of this formula defines the \emph{transitive closure} of the edge relation $E$. Consequently, the expression $\forall x_1\forall x_2\varphi^\infty(x_1,x_2)$ is a formula of least fixed-point logic LFP that expresses \emph{connectivity}.
 For a different example, let $\psi(x,S)$ be the formula $\forall y (E(y,x)\rightarrow T(y))$, where  $T$ is a unary relation symbol. It can be verified that, for every finite graph ${\bf G}=(V,E)$, the least fixed-point $\psi^\infty(x)$ defines the set of all nodes $v$ in $V$ such that no path containing $v$ leads to a cycle. Consequently, the expression $\forall x\psi^\infty(x)$ is a formula of least fixed-logic LFP that expresses \emph{acyclicity} on finite graphs.

Transitive Closure Logic (TCL) is the fragment of LFP that allows for the formation of the transitive closure of first-order definable relations. Thus, if $\theta(x_1,\ldots,x_k,x_{k+1},\ldots,x_{2k})$ is a first-order formula, then we can form in TCL the least fixed point of the formula:
\begin{multline*}\theta(x_1,\ldots,x_k,x_{k+1},\ldots,x_{2k}) \ \lor \\ \exists z_1\cdots\exists z_k(\theta(x_1,\ldots,x_k,z_1,\ldots,z_k)\land S(z_1,\ldots,z_k,x_{k+1},\ldots,x_{2k}))\end{multline*}

As regards their expressive power, it is known that FO $\subset$ TCL $\subset$ LFP on the class of all finite graphs.  In other words, FO is strictly less expressive than TLC, while TLC is strictly less expressive than LFP on the class of all finite graphs. As regards connections to computational complexity, it is known  FO is properly contained in LOGSPACE, TLC is properly contained in NLOGSPACE, and LFP is properly contained in PTIME on the class of all finite graphs.

The situation, however, changes if \emph{ordered} finite graphs are considered, that is, finite structures of the form ${\bf G}=(V,E,\leq)$, where $E$ is a binary relation on $V$ and $\leq $ is a total order on $V$ that can be used in LFP and in TCL formulas. In this setting, it is known that TLC = NLOGSPACE and that LFP = PTIME (the latter result is known as the Immerman-Vardi Theorem); thus, separating TLC from LFP on the class of all ordered graphs is equivalent to showing that NLOGSPACE is properly contained in PTIME, which is an outstanding open problem in computational complexity. Furthermore, similar results hold for the class of all finite structures and the class of all ordered finite structure over a relational schema containing at least one relation symbol of arity at least $2$. These results have been established in the context of \emph{descriptive complexity theory}, which studies the connections between computational complexity and expressibility in logics on finite structures. We refer the reader to the monographs \cite{DBLP:books/daglib/0095988,DBLP:books/sp/Libkin04} for detailed information.

\paragraph{Binary Non-dictatorial Aggregators and Transitive Closure Logic}
After the preceding digression into logic and complexity, we return to the question of when a domain $X$ of feasible evaluations
admits a binary non-dictatorial aggregation.

We begin by  first encoding a set $X$ of feasible evaluations by a suitable finite structure. To this effect,
we consider a relational schema $\bf R$ consisting of three unary relations $X'$, $I'$, $V'$, and one ternary relation $R'$. Intuitively, $X'$ will be interpreted by a set of feasible evaluations, $I'$ will be interpreted by the set of issues at hand, and $V'$ will be interpreted by  the set of positions over all issues.

Given a set $X\subseteq A^m$ of feasible evaluations, we let ${\bf A}(X) = (A,X,I,V)$ be the following finite $\bf R$-structure:
\begin{itemize}
\item   $A = X \cup I\cup V$, where $I= \{1,\ldots,m\}$ is the set of issues, and $V$ is the union of all positions over all issues.
    \item $R$ is the ternary relation consisting of all triples $(x,j,v)$ such that $x\in X$ and $v$ is the $j$-th coordinate of $x$, i.e., the position for issue $j$ in $x$.
    \end{itemize}

 It is clear that $X$ can be identified with the finite structure ${\bf A}(X)$. Conversely, if we are given a finite $\bf R$-structure $\bf A$ in which $R\subseteq X\times I\times V$, then $X$ can be thought of as a set of feasible evaluations over the issues $I$.
 \begin{lemma} \label{FO-lemma}
There is a first-order formula $\varphi(u,u',k,v,v',l)$ such that, for every set $X$ of feasible evaluations, we have that
  $\varphi(u,u',k,v,v',l)$ defines the edge relation of the directed graph $H_X$, when interpreted on the finite structure ${\bf A}(X)$.

\end{lemma}

\noindent\textbf{Proof}
Consider the following first-order formula
$\varphi(u,u',k,v,v',l)$:
\begin{align*}\exists & z\exists z' ((X'(z) \land X'(z') \land R'(z,k,u)\land R'(z,l,v)\land R(z',k,u')\land R(z',l,v')) \\\land &
\neg \exists y(X'(y) \land R'(y,k,u)\land R'(y,l,v')\\ \land & \forall j\forall w (R'(y,j,w) \rightarrow (R'(z,j,w) \lor R'(z',j,w))))).\end{align*}
It is immediate from the definition of  $H_X$ that the formula
 $\phi(u,u',k,v,v',l)$  defines indeed the edge relation of  $H_X$.
\hfill$\Box$\vspace{0.2cm}

We now have all the concepts and tools needed to obtain the following result.

\begin{theorem} \label{binary-TCL:cor}
The following problem is expressible in Transitive Closure Logic TCL: given a set $X$ of feasible evaluations (encoded as the finite structure ${\bf A}(X)$), does $X$ admit a binary, non-dictatorial aggregator?  Hence, this problem is also expressible in Least Fixed Point Logic LFP.
\end{theorem}

\noindent\textbf{Proof} The result follows immediately from Lemma \ref{lem:binary}, Lemma \ref{FO-lemma}, and the definition of Transitive Closure Logic.
\hfill$\Box$\vspace{0.2cm}

As stated earlier, every property that can be expressed in TCL is in \textrm{NLOGSPACE}. Thus, the problem of detecting if $X$ admits a binary non-dictatorial aggregator is in \textrm{NLOGSPACE}. Note that membership of this problem in \textrm{NLOGSPACE} could also be inferred from Lemma \ref{lem:binary} and the observation that the graph $H_X$ can be constructed in \textrm{LOGSPACE}. Lemma \ref{FO-lemma} strengthens this observation by showing that $H_X$ is actually definable in first-order logic, which is a small fragment of \textrm{LOGSPACE}.

Now, let $X\subseteq\{0,1\}^m$. We prove the following result. \begin{lemma}\label{lem:affinelogspce} Checking whether a domain $X\subseteq\{0,1\}^m$ is affine can be done in \textrm{LOGSPACE}. \end{lemma}
\textbf{Proof}
Suppose we have a Turing Machine with a read-only tape containing the tuples of $X$. In the work tape, we store triples $(i_1,i_2,i_3)$ of integers in $\{1,\ldots,n\}$ in binary. This takes $O(\log n)$ space. The integer $i_j$ points to the $i_j$-th element of $X$. We want to examine if the sum modulo 2 of these three elements is also in $X$.

To do that, for each such triple, we examine all integers $i_4 \leq n$ one at time. This adds another $\log n$ number of cells in the work tape. We then store all integers $j\leq m$ binary, one at a time, using another $\log m$ bits to the work tape.

Once we have $i_1$, $i_2$, $i_3$, $i_4$, and $j$ on the work tape, we check whether the entry $a^{i_4}_j=\oplus(a^{1_1}_j, a^{i_2}_j, a^{i_3}_j)$. If it is, we go to the next $j$. If it is not, we go to the next $i_4$, and when we are done with the triple $(i_1,i_2,i_3)$, we go to the next such triple. If for every triple $(i_1,i_2,i_3)$, we find a suitable integer $i_4$, $X$ is affine. Else, it is not.
\hfill$\Box$\vspace{0.2cm}

By Corollary \ref{cor:boolean}, Lemma \ref{lem:affinelogspce} and the discussion following Theorem \ref{binary-TCL:cor}, we obtain the following result.

\begin{theorem} \label{thm:boolean-tract} The following problem is in $\textrm{NLOGSPACE}$:  given a set $X \subseteq \{0,1\}^m$ of feasible evaluations in the Boolean framework, decide whether or not $X$ is a possibility domain.
\end{theorem}

\subsection{Tractability of Uniform Possibility Domains}\label{ssec:upd-tract}
Recall that a constraint language is  a finite set $\Gamma$ of relations of finite arities over a finite non-empty set $A$. The {\em conservative constraint satisfaction problem for} $\Gamma$, denoted by $\mbox{c-CSP}(\Gamma)$, is the constraint satisfaction problem for the constraint language $\overline{\Gamma}$ that consists of the relations in $\Gamma$ and, in addition, all unary relations on $A$. Intuitively, this amounts to the ability to arbitrarily restrict the domain of each variable in a given instance.

Bulatov \cite{bulatov2006dichotomy,bulatov2011complexity} established a dichotomy theorem for the computational complexity of $\mbox{c-CSP}(\Gamma)$: if for every two-element subset $B$ of $A$, there is a conservative polymorphism $f$ of $\Gamma$ such that  $f$ is binary and $f{\restriction_B}\in \{\wedge, \vee\}$ or $f$ is ternary and $f{\restriction_B}\in \{{\rm  maj}, \oplus\}$, then $\mbox{c-CSP}(\Gamma)$ is solvable in polynomial time; otherwise, $\mbox{c-CSP}(\Gamma)$ is NP-complete.
Carbonnel showed that the boundary of the dichotomy for $\mbox{c-CSP}(\Gamma)$ can be checked in polynomial time.
\begin{customthm}{E}{\cite[Theorem 4]{carbonnel2016dichotomy}}\label{thm:carb-tract}
There is a polynomial-time algorithm for the following problem: given a constraint language $\Gamma$ on a set $A$, determine whether or not for every two-element subset $B\subseteq A$, there is a conservative polymorphism $f$ of $\Gamma$ such that either $f$ is binary and $f{\restriction B}\in \{\wedge, \vee\}$ or
$f$ is ternary and $f{\restriction B} \in \{{\rm  maj}, \oplus\}$. Moreover, if such a polymorphism exists, then the algorithm produces one in polynomial time.
\end{customthm}
The final results of this section is about the complexity of detecting uniform possibility domains.
\begin{theorem}\label{thm:tractupd}
There is a polynomial-time algorithm for solving the following problem:  given a set  $X$  of feasible evaluations, determine whether or not $X$ is a {\un} possibility domain and, if it is,  produce  a ternary weak near-unanimity aggregator for $X$.
\end{theorem}
In what follows, given a two-element set $B$, we will arbitrarily identify its elements with $0$ and $1$.
Consider the functions $\wedge^{(3)}$ and $\vee^3$ on $\{0,1\}^3$, where $\wedge^3(x,y,z):=(\wedge(\wedge(x,y),z))$ and $\vee^3(x,y,z):=(\vee(\vee(x,y),z))$. It is easy to see that the only ternary, conservative, weak near-unanimity functions on $\{0,1\}$ are $\wedge^{(3)}$, $\vee^3$, ${\rm maj}$, and $\oplus$. We will also make use of the following lemma, which gives an alternative formulation  of the boundary of the dichotomy for conservative  constraint satisfaction.

\begin{lemma}\label{lemma:wnu}
Let $\Gamma$ be a constraint language on set $A$. The following two statements are equivalent.
\begin{enumerate}
    \item For every two-element subset $B\subseteq A$, there exists a conservative polymorphism $f$ of $\Gamma$ (which, in general, depends on $B$), such that $f$ is binary and $f{\restriction_B}\in\{\wedge,\vee\}$  or $f$ is ternary and $f{\restriction_B}\in\{{\rm maj},\oplus\}$.\label{cond:one}
    \item $\Gamma$  has a ternary, conservative, weak near-unanimity polymorphism.
\end{enumerate}
\end{lemma}

\noindent\textbf{Proof} (\emph{Sketch})
($1\Rightarrow 2$)
Given a two-element subset $B\subseteq A$ and a binary conservative  polymorphism $f$ of $\Gamma$ such that $f{\restriction_B}\in\{\wedge,\vee\}$, define $f'$ to be the ternary operation such that $f'(x,y,z)=f(f(x,y),z)$,  for all $x,y,z\in A$.
It is easy to see that $f'$ is a conservative polymorphism of $\Gamma$ as well and also that
 $f'{\restriction_B}\in\{\wedge^{(3)},\vee^{(3)}\}$.

The hypothesis and the preceding argument imply that, for each two-element subset $B\subseteq A$, there exists a ternary conservative  polymorphism $f$ of $\Gamma$ (which, in general,  depends on $B$)  such that $f{\restriction_B}\in\{\wedge^{(3)},\vee^{(3)},{\rm maj},\oplus\}$. For each two-element subset $B\subseteq A$, select such a polymorphism and let $f^1,\ldots,f^N$, $N\geq1$, be an enumeration of all these polymorphisms. Clearly, the restriction of each $f^i$ to its respective two element subset is a weak near-unanimity operation.

Consider the `$\diamond$' operator that takes as input  two ternary operations $f,g:A^3\to A$  and returns as output a ternary operation $f\diamond g$ defined by $$(f\diamond g)(x,y,z):=f(g(x,y,z),g(y,z,x),g(z,x,y)).$$
If $f,g$ are conservative polymorphisms of $\Gamma$, then so is $(f\diamond g)$. Also, if $B$ is a two-element subset of $A$ such that   $f{\restriction_B}$ or $g{\restriction_B}$  is a weak near-unanimity operation, then so is
$(f\diamond g){\restriction_B}$. Consider now the iterated diamond operation $h$ with
$$h:=f^1\diamond(f^2\diamond(\ldots\diamond(f^{N-1}\diamond f^N)\ldots)).$$ By the preceding discussion,  $h$ is a conservative polymorphism such that $h{\restriction_B}$ is a weak near-unanimous operation for \emph{every}  two-element subset $B$ of $A$, hence  $h$ itself  is a weak near-unanimity, conservative, ternary operation of $\Gamma$.

($2\Rightarrow 1$) Let $h$ be a ternary, conservative, weak near-unanimity polymorphism of $\Gamma$. Thus, for every two-element subset $B\subseteq A$, we have that $h{\restriction_B}\in\{\wedge^{(3)},\vee^{(3)},{\rm maj},\oplus\}.$

If there is a two-element subset $B\subseteq A$ such that $h{\restriction_B}\in\{\wedge^{(3)},\vee^{(3)}\}$, then consider the binary function $g$, which is defined as: $$g(x,y):=h(x,x,y)=h(pr_1^2(x,y),pr_1^2(x,y),pr_2^2(x,y)).$$
Obviously, $g$ is a binary conservative polymorphism of $\Gamma$; moreover, for every two-element subset $B\subseteq A$, if
  $h{\restriction_B}\in\{\wedge^{(3)},\vee^{(3)}\}$, then $g{\restriction_B}\in\{\wedge,\vee\}$. \hfill$\Box$\vspace{0.2cm}

For a detailed proof of Lemma \ref{lemma:wnu}, see \cite[Theorem 5.5]{kirousis2019aggregation}.\vspace{0.2cm}

\noindent\textbf{Proof of Theorem \ref{thm:tractupd}} By Theorem \ref{thm:updcar}, a set $X$ of feasible evaluations is a {\un} possibility domain if and only if
there is a ternary aggregator $F=(f_1,\ldots,f_m)$ such that each $f_j$ is a weak near-unanimity operation, i.e., for all $j \in \m$  and for all $x,y \in X_j$, we have that
$f_j(x,y,y) = f_j(y,x,y)= f_j(y,y,x)$.
As in the proof of Theorem \ref{thm:pdmain}, we can go back-and-forth between $X$ and the set $\tilde{X}$ and verify that $X$ is a {\un} possibility domain if and only if $\tilde{X}$ has a ternary, conservative, weak near-unanimity polymorphism. Theorem \ref{thm:carb-tract} and Lemma \ref{lemma:wnu} then imply that the existence of such a polymorphism can be tested in polynomial time, and that such a polymorphism  can be produced in polynomial time, if one exists.\hfill$\Box$\vspace{0.2cm}

In the Boolean case, we can prove the tractability of detecting locally non-dictatorial aggregators without using Theorem \ref{thm:carb-tract}. This will allow us to obtain better complexity bounds in Section \ref{sec:implicit}. In the Boolean case, Theorem \ref{thm:updcar} has been strengthened by Diaz et al.: \begin{corollary}{\cite[Corollary $4.1$]{diaz2019syntactic}}\label{cor:lpdchar}
A Boolean domain $X\subseteq\{0,1\}^m$ is a local possibility domain if and only if it admits a ternary aggregator $F=(f_1,\ldots,f_m)$ such that $f_j\in\{\wedge^{(3)},\vee^{(3)},\oplus\}$, for $j=1,\ldots,m$.
\end{corollary}
Thus, we can obtain the following algorithm in the Boolean case.
\begin{corollary}\label{cor:tractlpd}
There is a polynomial-time algorithm for solving the following problem: given a Boolean set $X\subseteq\{0,1\}^m$ of feasible evaluations, determine whether or not $X$ is a local possibility domain and, if it is,  produce  a ternary aggregator $F=(f_1,\ldots,f_m)$ for $X$ such that $f_j\in\{\wedge^{(3)},\vee^{(3)},\oplus\}$, $j=1,\ldots,m$.
\end{corollary}
\textbf{Proof} Let $F=(f_1,\ldots,f_m)$ be the binary maximum symmetric aggregator obtained by Corollary \ref{cor:maxsym} and let $I,J\subseteq\m$ such that $f_i=\wedge$ for all $i\in I$ and $f_j=\vee$ for all $j\in J$ (both $I$ and $J$ can be empty). We prove that $X$ is a local possibility domain if and only if it admits the ternary aggregator $G=(g_1,\ldots,g_m)$, where $g_i=\wedge^{(3)}$ for all $i\in I$, $g_j=\vee^{(3)}$ for all $j\in J$ and $g_k=\oplus$, for all $k\in\m\setminus(I\cup J)$. Since we can obtain $F$ in polynomial-time in the size of $X$ and since checking whether $G$ is an aggregator for $X$ can be done also in polynomial time, the procedure is clearly polynomial in the size of $X$.

That $X$ is a local possibility domain if it admits $G$ is self-evident. Assume now that $X$ is a local possibility domain. Let also $F^{(3)}=(f^{(3)}_1,\ldots,f^{(3)}_m)$ be the $m$-tuple of ternary functions, where $$f_j^{(3)}(x,y,z)=f_j(f_j(x,y),z),$$ $j=1,\ldots,m$. It is straightforward to check that (i) $F^{(3)}$ is an aggregator for $X$, (ii) $f_i^{(3)}=\wedge^{(3)}$ for all $i\in I$, (iii) $f_j^{(3)}=\vee^{(3)}$ for all $j\in J$ and (iv) $f_k^{(3)}\in\{\pr_1^3,\pr_3^3\}$ for all $k\in\m\setminus(I\cup J)$.

Now, since $X$ is a local possibility domain, by Corollary \ref{cor:lpdchar}, $X$ admits a ternary aggregator $H=(h_1,\ldots,h_m)$, such that $h_j\in\{\wedge^{(3)},\vee^{(3)},\oplus\}$, $j=1,\ldots,m$. If $H=G$, there is nothing to prove.

First, assume that there is some $k\notin I\cup J$, such that $h_k\in\{\wedge^{(3)},\vee^{(3)}\}$. Now, let $G'=(g'_1,\ldots,g'_m):=H\diamond F^{(3)}$ and $F'=(f'_1,\ldots,f'_m)$ be a binary $m$-tuple of functions such that: $$f'_j(x,y):=g'(x,x,y),\text{ for all }x,y\in A.$$ It holds that $G'$ is an aggregator for $X$ such that $g'_j\in\{\wedge^{(3)},\vee^{(3)}\}$, for all $j\in I\cup J\cup\{k\}$. Thus, $F'$ is a non-dictatorial aggregator and $f'_j$ is symmetric for all $j\in I\cup J\cup\{k\}$. Contradiction, since $F$ is a maximum symmetric aggregator and $f_k$ is not symmetric.

Finally, suppose that there is some $j\in I\cup J$ such that $h_k=\oplus$, for all $k\in (\{1,\ldots,m\}\setminus(I\cup J))\cup\{j\}$. Then, $G=H\diamond F^{(3)}$ and thus is an aggregator for $X$.

\section{Implicitly given Domains}\label{sec:implicit}
In Subsection \ref{ssec:logicbased}, we describe the logic-based approach, where the domain $X$ is given implicitly via an integrity constraint or an agenda. In Subsections \ref{ssec:integr} and \ref{ssec:agenda}, we establish complexity bounds for checking whether a domain is a  possibility domain or a local possibility domain in  these two variants of the logic-based approach (recall that a local possibility domain is a uniform possibility domain in the Boolean framework). Finally, in Subsection \ref{ssec:other}, we extend these results to other types of non-dictatorial aggregation that have been used in the literature.

\subsection{The Logic-based Approach}\label{ssec:logicbased}
In this subsection, we present two ways that a set $X$ of feasible voting patterns can be given implicitly: as an integrity constraint and as an agenda. Both variants are in the Boolean framework and have been studied extensively in the literature.

Suppose that we have a \emph{propositional formula} $\phi$ on $m$ variables $x_1,\ldots,x_m$. Each variable $x_j$ corresponds to the $j$-th issue, $j=1,\ldots,m$, where the possible positions are $0$ and $1$. Let $X_{\phi}:=\rMod(\phi)$ be the set consisting of all $m$-ary vectors of \emph{satisfying truth assignments} of $\phi$. In this setting, we say that $\phi$ is an \emph{integrity constraint}.

For the second variant, suppose that we have an \emph{agenda} $\ph=(\phi_1,\ldots,\phi_m)$ of $m$ propositional formulas. For a formula $\psi$ and an $x\in\{0,1\}$, let
$$\psi^x := \begin{cases} \psi & \mbox{ if } x=1,\\
\neg\psi & \mbox{ if } x=0.
\end{cases}.$$
Finally, let:
$$X_{\ph} := \left\{\x= (x_1, \ldots, x_m) \in \{0,1\}^m \mid \bigwedge_{j=1}^m {\phi}_j^{x_j} \mbox{ is satisfiable } \right\}.$$

Recall that, as usual in aggregation theory, we have assumed that domains $X$ are non-degenerate, i.e., $|X_j|\geq 2$ (thus $X_j=\{0,1\}$ in the Boolean framework), for $j=1,\ldots,m$. Thus, we assume that both integrity constraints and agendas are such that their domains are non-degenerate. On the other hand, when we consider (propositional) formulas, we do not assume anything regarding their domain (it can even be empty).

It is well known that given a domain $X\subseteq\{0,1\}^m$, there is a formula $\phi$ such that its set of models $\rMod(\phi)$ is equal to $X$; see, e.g., \cite{enderton2001mathematical}. Dokow and Holzman \cite{dokow2010aggregationnonB} prove that there is also an agenda $\ph$ such that $X_{\ph}=X$. Thus, the three variants (explicit representation, implicit representation via an integrity constraint, and implicit representation via an agenda) are in some sense equivalent, as regards the existence of (local) non-dictatorial aggregators. However, neither the integrity constraint nor the agenda describing a given domain need  be unique. Thus, there is a possible loss of information when passing from one variant to another; as List and Puppe \cite{list2009judgment} argue, this can be significant for certain aspects of the aggregation problem. Consider, for example, the agendas $$\bar{\phi_1}=\{p,q,p\wedge q\}$$ and $$\bar{\phi_2}=\{p\vee q\vee\neg r, p\vee\neg q\vee\neg r,\neg p\vee q\vee\neg r,\neg p\vee\neg q\vee r\}.$$
It is easy to see that $$X_{\ph_1}=X_{\ph_2}=\{(0,0,0),(0,1,0),(1,0,0),(1,1,1)\}.$$ Nevertheless, the two agendas can have different meaning regarding the problems they model. First of all, the formulas of $\ph_1$ consist  of two independent propositional variables, while those of $\ph_2$ by three. Furthermore, the fact that $\ph_1$ contains propositional variables, whereas $\ph_2$ does not, can also lead to different strategies in order to use these agendas, as in \cite{mongin2008factoring}, where Mongin replaces IIA with the \emph{Independence of Irrelevant Propositional Alternatives} axiom.

In terms of the computational complexity of passing from one framework to another, Zanuttini and H\'ebrard \cite{zanuttini2002unified} show that given a domain $X$, one can construct a formula $\phi$ such that $\rMod(\phi)=X$ in polynomial time in the size of the domain. Also, the construction in \cite{dokow2009aggregation}, where given a domain $X$, we obtain an agenda $\bar{\phi}$ such that $X_{\bar{\phi}}=X$ can obviously be carried out in polynomial time in the size of the domain. It is very easy to find integrity constraints and agendas whose domains are exponentially large on their respective sizes: consider for example the integrity constraint $(x_1\vee\neg x_1)\wedge\cdots\wedge(x_m\vee\neg x_m)$ and the agenda $(x_1,\ldots,x_m),$ where $x_1,\ldots,x_m$ are pairwise distinct variables. Both have domains equal to the full Boolean domain $\{0,1\}^m$. Finally, Endriss et al. \cite[Proposition 9]{endriss2016succinctness} show that, unless the \emph{polynomial hierarchy} collapses, we cannot describe any arbitrary agenda by an integrity constraint of polynomial size to that of the agenda and that, given an integrity constraint $\phi$, the problem of finding an agenda $\ph$ such that $X_{\ph}=X_{\phi}$ is $\FNP$-complete.

Here, we examine the computational complexity of checking if a domain $X$ is a (local) possibility domain in both the integrity constraint variant and the agenda variant. In all that follows, we assume that the integrity constraints are defined on at least three variables and agendas contain at least three propositional formulas, since domains $X\subseteq\{0,1\}^m$ where $m=1$ or $2$ are all possibility domains.

\subsection{Integrity Constraints}\label{ssec:integr}
Let $\phi$ be an integrity constraint on $m$ variables and let $X_{\phi}=\rMod(\phi) \subset \{0,1\}^m$. The following theorems provide upper and lower bounds to the complexity of checking if $X_{\phi}$ is a (local) possibility domain. For the upper bounds, we work with oracles, using the definition of the graph $H_X$ in Section \ref{sec:abstract} and make use of the following straightforward fact.

\begin{lemma}\label{lem:closure}
Let $F=(f_1,\ldots,f_m)$ be an $m$-tuple of $n$-ary conservative and polynomial-time computable functions. Deciding, on input $\phi$, whether $F$ is an aggregator for $X_{\phi}$ is in $\rP_1^P=\coNP$.	
\end{lemma}
\textbf{Proof} The result is immediate since the problem can be cast as follows:
\begin{quote}
    for all $m$-tuples $\x^1,\ldots,\x^n\in\{0,1\}^m$, if all $n$ satisfy $\phi$, then so does the $m$-tuple $F(\x^1,\ldots,\x^n)$
\end{quote}
and since checking if a formula $\phi$ is satisfied by a specific assignment can be done in polynomial time.\hfill$\Box$\vspace{0.2cm}

A function $f$ is polynomial-time computable if given its input $\bar{x}$, we can compute its output $f(\bar{x})$ in polynomial time. For our purposes, it suffices that the functions $\wedge$, $\vee$, $\oplus$, $\wedge^{(3)}$, $\vee^{(3)}$ and $\pr_i^n$ are all polynomial-time computable, for all $n\in\N$ and $i\in\n$.

In terms of lower bounds, we provide polynomial-time reductions from two $\coNP$ -complete problems: the \emph{semantical independence} problem and the \emph{unsatisfiability} problem for propositional formulas. The latter is the well known problem of whether a formula has no satisfying assignments. The former asks whether a given propositional formula is \emph{(semantically) dependent} to all its variables; that is, there is no variable (or set of variables) such that whether an assignment of values satisfies the given formula or not, does not depend on the values of the variable(s). For a systematic overview of this problem and its variations, see \cite{lang2002conditional} and \cite{lang2003propositional}. In the setting of agendas, this notion has been studied under the name \emph{agenda separability} (see \cite{lang2016agenda}).

For domain $X\subseteq\{0,1\}^m$ and a nonempty subset $I\subset\{1,\ldots,m\}$, we denote by $X_I$ the projection of $X$ to $I$, that is the set of all partial vectors with indices in $I$ that can be extended to elements of $X$. Let also $X_{-I}:=X_{\m\setminus I}$ and $X\approx Y$ mean that we can obtain $X$ by permuting the columns of $Y$.
\begin{definition}\label{def:indep}
Let $\phi(x_1,\ldots,x_m)$ be a propositional formula, where $X:=\rMod(\phi)$ and let $V\subseteq\{x_1,\ldots,x_m\}$ be a subset of its variables. Let also $i\in\m$. We say that $\phi$ is:\begin{itemize}
    \item[i.] (semantically) independent from variable $x_i$ if: $$X\approx X_{\{i\}}\times X_{-\{i\}}$$
    \item[ii.](semantically) independent from the set of variables $V$ if it is independent from every $x_j\in V$.
\end{itemize}
\end{definition}
In our setting, an integrity constraint being independent from a variable $x_j$ means that issue $j$ does not contribute anything in the logical consistency restrictions imposed by the constraint. Lang et al. \cite{lang2003propositional} showed that the problem of checking if a propositional formula depends on all its variables (is \emph{simplified variable-dependent} in their terminology), is $\coNP$-complete.

To make our reductions easier to follow, we work with the specific domain: $$\imp:=\{0,1\}^3\setminus\{(0,0,0),(1,1,1)\}.$$
Observe that $\imp$ corresponds to a natural and well studied problem in both preference and judgment aggregation. Suppose that we have three alternatives $A$, $B$ and $C$, where issue $1$ corresponds to deciding between $A$ and $B$, issue $2$ corresponds to deciding between $B$ and $C$, and issue $3$ corresponds to deciding between $C$ and $A$. In that setting, we can assume that, in each issue, $1$ denotes preferring the former option and $0$ the latter. One can easily see now that $\imp$ corresponds to the natural requirement of transitivity of preferences.

Consider now the following lemma, which is a version of Arrow's impossibility result in Judgment Aggregation.
\begin{lemma}\label{lem:imp}
$\imp$ is an impossibility domain.
\end{lemma}
\textbf{Proof} By Corollary \ref{cor:boolean}, we only need to check if $\imp$ is affine, or if it admits a binary non-dictatorial aggregator. Easily. $$\bar{\oplus}((1,0,0),(0,1,0),(0,0,1))=(0,0,0)\notin \imp,$$ thus $\imp$ is not affine.

On the other hand, let $F=(f_1,f_2,f_3)$ be a binary non-dictatorial triple of functions. There are $4^3-2=62$ cases for $F$. We arbitrarily choose to show three of them. The rest are left to the interested reader.\begin{itemize}
    \item If $f_1=f_2=\wedge$ and $f_3=\vee$, then $F((1,0,0),(0,1,0))=(0,0,0)\notin \imp$.
    \item If $f_1=\wedge$, $f_2=\vee$ and $f_3=\pr_1^2$, then $F((1,0,0),(0,0,1))=(0,0,0)\notin \imp$.
    \item If $f_1=\vee$, $f_2=\pr_1^2$ and $f_3=\pr_2^2$, then $F((0,0,1),(0,1,0))=(0,0,0)\notin \imp$.
\end{itemize}
Thus, for $F$ to be an aggregator for $\imp$, it must hold that $f_1=f_2=f_3=\pr_d^2$, $d=1,2$ and, consequently $\imp$ admits only dictatorial binary aggregators. \hfill$\Box$\vspace{0.2cm}

Let $\psi$ be the propositional formula:\begin{equation}\label{impic}
    \psi(y_1,y_2,y_3)=(y_1\vee y_2\vee y_3)\wedge(\neg y_1\vee\neg y_2\vee\neg y_3).
\end{equation}
Easily, $\rMod(\psi)=\imp$. We are now ready to obtain our results.

\begin{theorem}\label{thm:possintcon}
Deciding, on input $\phi$,  whether $X_{\phi}$  admits a non-dictatorial aggregator is (i) in $\rS^P_2\cap\rP_2^P$ and (ii) $\coNP\text{-hard}$.	
\end{theorem}
\textbf{Proof}
(i) By Corollary \ref{cor:boolean}, $X_{\phi}$ is a possibility domain if and only if it admits a binary non-dictatorial aggregator or it is affine. The problem of whether $X_{\phi}$ is affine is in $\rP_1^P=\coNP$ (and thus in $\rS_2^P\cap\rP_2^P$ too), by Lemma \ref{lem:closure}, using $\bar{\oplus}$ for the given $m$-tuple of functions (recall also Lemma \ref{sch:affine}).

For the problem of the existence of binary non-dictatorial aggregators for $X_{\phi}$, we show separately that it is both in $\rS_2^P$ and in $\rP_2^P$. For the former, note that there are only four conservative (equivalently Paretian) functions from $\{0,1\}^2\mapsto\{0,1\}$, namely $\pr_1^2$, $\pr_2^2$, $\wedge$ and $\vee$. Therefore, there are $4^m-2$ tuples   $F = (f_1, \ldots f_m)$, with $f_j: \{0,1\}^2 \mapsto \{0,1\}, j=1, \ldots m$ of $m$  conservative functions, where not all $f_j$ are the same projection function. Such an $m$-tuple can be thought of as a binary $(2\times 3)m$-sequence, where each $f_j$ is encoded by a sequence $(01a10b)$, $a,b\in\{0,1\}$, meaning that $f_j(0,1)=a$ and $f_j(1,0)=b$.

Let us call such $F$ {\em candidates for non-dictatorial aggregators}.  The question of deciding whether a given   $F$ is one of the $4^m -2$  candidates for non-dictatorial aggregators is easily in $\rm{P}$. Also, the question of whether a given binary $F$ is an aggregator for $X_{\phi}$ is again in $\rP_1^P$ by Lemma \ref{lem:closure}. Therefore the problem of whether $X_{\phi}$ admits a binary non-dictatorial aggregator is in $\rm{\Sigma}^P_2$ because it can be cast as:
\begin{quote}
There exists  a $F = (f_1, \ldots f_m)$ such that $F $ is a candidate for non-dictatorial aggregator and for all $m$-tuples $\x, \y \in \{0,1\}^m$, if both satisfy $\phi$, then so does the $m$-tuple $F(\x,\y)$.
\end{quote}

To show that it is also in $\rP_2^P$, recall that, by Lemma \ref{lem:binary}, the set $X_{\phi}$ admits a binary non-dictatorial aggregator if and only if the graph $H_X$ is not strongly connected. We will show that checking if $H_{X_{\phi}}$ is strongly connected is in $\rS_2^P$, which means that checking if $H_{X_{\phi}}$ is not strongly connected is in $\rP_2^P$.

First note that the size of $H_X$ is polynomial in the size of $\phi$, since it has $2m$ nodes, where $m$ is the number of variables of $\phi$. Thus, it suffices to prove that testing whether two nodes of $H_{X_{\phi}}$ are connected is in $\NP$ with an oracle in $\coNP$.

To test if two nodes are connected, we first obtain a path witnessing this. To verify it is indeed a path, we need to check if any two of its consecutive nodes, say $uu'_s$ and $vv'_t$, are connected by the edge $uu'_s\rightarrow vv'_t$ in $H_X$. To do that, we can again take the satisfying assignments $z$ and $z'$ of $\phi$ that witness that $uu'_s\rightarrow vv'_t$. Then, using the $\coNP$ oracle, we need to check that there is no $z^*$ that: (i) satisfies $\phi$, (ii) extends $uv'$ and (iii) agrees on every coordinate either $z$ or $z'$.

(ii) Let $\chi(x_1,\ldots,x_k)$ be a propositional formula on $k$ variables. We construct, in polynomial time, a formula $\phi$ such that $\chi$ is independent from at least one of its variables if and only if $X_{\phi}$ is a possibility domain. In all that follows, $m=k+3$.

Let: $$\phi(x_1,\ldots,x_k,y_1,y_2,y_3)=
            \chi(x_1,\ldots,x_l)\oplus\psi(y_1,y_2,y_3),$$ where: $\{x_1,\ldots,x_k\}$ and $\{y_1,y_2,y_3\}$ are disjoint sets of variables.

First, note that the length $|\phi|$ of $\phi$ is linear to that of $\chi$, since $|\phi|=|\chi|+6$. Thus the construction is polynomial.

By \eqref{impic}, it holds that: \begin{equation}\label{modx}
X_{\phi}=\Big(\rMod(\chi)\times\{(0,0,0),(1,1,1)\}\Big)\cup\Big(\rMod(\neg\chi)\times\imp\Big).\end{equation}

\begin{claim}\label{claim:indvar}
$\chi$ is independent from at least one of its variables if and only if $X_{\phi}$ is a possibility domain.
\end{claim}
\textit{Proof of Claim:}
We first consider the two extreme cases. If $\chi$ is unsatisfiable or a tautology, then by \eqref{modx} we have that:$$X_{\phi}=\{0,1\}^k\times\imp,$$ or that $$X_{\phi}=\{0,1\}^k\times\{(0,0,0),(1,1,1)\}$$ respectively. In both cases, we have that $\chi$ is independent from all its variables and $X_{\phi}$ is a possibility domain, since it is a Cartesian product (recall Lemma \ref{lem:cartprod}).

Thus, we can assume that both $\rMod(\chi)$ and $\rMod(\neg\chi)$ are not empty. We proceed with a series of claims.
\begin{claim}\label{claim:affine}
$X_{\phi}$ is not affine.
\end{claim}
\textit{Proof of Claim:} Let $\bar{a}:=(a_1,\ldots,a_k)\in\rMod(\neg\chi)$. Then, $(\bar{a},0,1,1)$, $(\bar{a},1,0,1)$, $(\bar{a},1,1,0)\in X_{\phi}$ and: $$\bar{\oplus}((\bar{a},0,1,1),(\bar{a},1,0,1),(\bar{a},1,1,0))=(\bar{a},0,0,0)\notin X_{\phi}.$$ Thus, $X_{\phi}$ is not affine. \hfill$\Box$\vspace{0.2cm}

By Corollary \ref{cor:boolean} and Claim \ref{claim:affine}, $X_{\phi}$ is a possibility domain if and only if it admits a binary non-dictatorial aggregator. The following claims show that such an aggregator must be in a restricted class.
\begin{claim}\label{claim:imp}
Assume $F=(f_1,\ldots,f_m)$ is a binary aggregator for $X_{\phi}$. Then, $f_{k+1}=f_{k+2}=f_{k+3}=\pr_d^2$, $d\in\{1,2\}$.
\end{claim}
\textit{Proof of Claim:} To obtain a contradiction, assume $(f_{k+1},f_{k+2},f_{k+3})\neq(\pr_d^2,\pr_d^2,\pr_d^2)$, $d=1,2$. By Lemma \ref{lem:imp}, $\imp$ is an impossibility domain. Thus, there exist $\bar{x},\bar{y}\in\imp$ such that $$\bar{z}:=(f_{k+1},f_{k+2},f_{k+3})(\bar{x},\bar{y})\in\{(0,0,0),(1,1,1)\}.$$ Let $\bar{a}:=(a_1,\ldots,a_k)\in\rMod(\neg\chi)$. Then $(\bar{a},\bar{x})$, $(\bar{a},\bar{y})\in X_{\phi}$, but: $$F((\bar{a},\bar{x}),(\bar{a},\bar{y}))=(\bar{a},\bar{z})\notin X_{\phi}.$$ Thus, $F$ is not an aggregator for $X_{\phi}$. Contradiction. \hfill$\Box$ \vspace{0.2cm}

The following claim states that $F=(f_1,\ldots,f_k,f_{k+1},f_{k+2},f_{k+3})$ cannot have its first $k$ coordinates be projections to the same coordinate $d\in\{1,2\}$ and the last three be projections to the other. This can also be derived by \eqref{modx}, since $X_{\phi}$ is not a Cartesian product. For a proof of this general and straightforward characterization, the interested reader is referred to \cite{diaz2019syntactic}. Here we opted to showcase the technique we follow throughout the rest of the proof.
\begin{claim}\label{claim:proj}
Assume $F=(f_1,\ldots,f_m)$ is a binary aggregator for $X_{\phi}$. If $f_{k+1}=f_{k+2}=f_{k+3}=\pr_d^2$, then there is at least one $j\in\{1,\ldots,k\}$ such that $f_j\neq\pr_{d'}^2$, $d,d'\in\{1,2\}$, $d\neq d'$.
\end{claim}
\textit{Proof of Claim:} We show the claim for $d=2$ and $d'=1$. The analogous arguments hold for the case where $d=1$ and $d'=2$.

To obtain a contradiction, assume that $f_1=\cdots=f_k=\pr_1^2$. Let $\bar{a}:=(a_1,\ldots,a_k)\in\rMod(\chi)$ and $\bar{b}:=(b_1,\ldots,b_k)\in\rMod(\neg\chi)$. Then, $(\bar{a},0,0,0)$ and $(\bar{b},0,0,1)\in X_{\phi}$, but: $$F((\bar{a},0,0,0),(\bar{b},0,0,1))=(\bar{a},0,0,1)\notin X_{\phi}.$$ Thus, $F$ is not an aggregator for $X_{\phi}$. Contradiction. \hfill$\Box$ \vspace{0.2cm}

\begin{claim}\label{claim:sym}
Assume $F=(f_1,\ldots,f_m)$ is a binary aggregator for $X_{\phi}$. Then, there is at least one $j\in\{1,\ldots,k\}$ such that $f_j$ is not symmetric.
\end{claim}
\textit{Proof of Claim:} By Claim \ref{claim:imp}, for some $d\in\{1,2\}$, $f_{k+1}=f_{k+2}=f_{k+3}=\pr_d^2$. To obtain a contradiction, assume that $f_j$ is symmetric for all $j\in\{1,\ldots,k\}$. Assume also w.l.o.g. that $d=2$. The analogous arguments work for $d=1$.

Let $\bar{a}:=(a_1,\ldots,a_k)\in\rMod(\chi)$, $\bar{b}:=(b_1,\ldots,b_k)\in\rMod(\neg\chi)$ and $$(f_1,\ldots,f_k)(\bar{a},\bar{b}):=\bar{c}.$$ Then, $(\bar{a},0,0,0)$ and $(\bar{b},0,0,1)\in X$. Since $F$ is an aggregator for $X_{\phi}$: \begin{align*}
    f((\bar{a},0,0,0),(\bar{b},0,0,1)) & =(\bar{c},0,0,0)\in X_{\phi},\\
    f((\bar{b},0,0,1),(\bar{a},0,0,0)) & =(\bar{c},0,0,1)\in X_{\phi},
\end{align*} which, by \eqref{modx}, implies that $\bar{c}\in\rMod(\chi)\cap\rMod(\neg\chi)$. Contradiction. \hfill$\Box$ \vspace{0.2cm}

The last claim deals with the case where we have both symmetric and non-symmetric components in $(f_1,\ldots,f_k)$. It completely outlines the class of binary aggregators available for $X_{\phi}$. Notationally, if $\bar{a}\in\{0,1\}^m$ and $I\subseteq\m$, $\bar{a}_I$ denotes the projection of $\bar{a}$ to the coordinates in $I$.
\begin{claim}\label{claim:class}
Assume $F=(f_1,\ldots,f_m)$ is a binary aggregator for $X_{\phi}$. Then, there exists a non-empty subset $J\subseteq\{1,\ldots,k\}$ such that $f_j=\pr_d^2$ for all $j\in J\cup\{k+1,k+2,k+3\}$, $d=1,2$.
\end{claim}
\textit{Proof of Claim:} By Claim \ref{claim:imp}, for some $d\in\{1,2\}$, $f_{k+1}=f_{k+2}=f_{k+3}=\pr_d^2$. To obtain a contradiction, assume that $f_j\neq\pr_d^2$, for all $j\in\{1,\ldots,k\}$. Assume also w.l.o.g. that $d=2$. The analogous arguments work for $d=1$.

By Claims \ref{claim:proj} and \ref{claim:sym}, there exists a partition $(I,J)$ of $\{1,\ldots,k\}$, such that $f_i$ is symmetric for all $i\in I$ and $f_j=\pr_1^2$ for all $j\in J$. To make things easier to follow, assume w.l.o.g. that there exists an $s\in\{1,\ldots,k-1\}$ such that $I=\{1,\ldots,s\}$ and $J=\{s+1,\ldots,k\}$. Let $\bar{a}:=(a_1,\ldots,a_k)\in\rMod(\chi)$, $\bar{b}:=(b_1,\ldots,b_k)\in\rMod(\neg\chi)$ and assume that: $$(f_1,\ldots,f_s)(\bar{a}_I,\bar{b}_I):=\bar{c}.$$ Then, $(\bar{a},0,0,0)$ and $(\bar{b},0,0,1)\in X_{\phi}$. Since $F$ is an aggregator for $X_{\phi}$, it must hold that: $$F((\bar{a},0,0,0),(\bar{b},0,0,1))=(\bar{c},\bar{a}_J,0,0,1)\in X_{\phi},$$ which, by \eqref{modx}, implies that $(\bar{c},\bar{a}_J)\in\rMod(\neg\chi)$. Furthermore, again since $F$ is an aggregator for $X_{\phi}$, it must be the case that: $$F((\bar{c},\bar{a}_J,0,0,1),(\bar{a},0,0,0))=(\bar{c},\bar{a}_J,0,0,0)\in X_{\phi},$$ which, by \eqref{modx}, implies that $(\bar{c},\bar{a}_J)\in\rMod(\chi)$. Thus, $(\bar{c},\bar{a}_J)\in\rMod(\chi)\cap\rMod(\neg\chi)$. Contradiction. \hfill$\Box$\vspace{0.2cm}

By Claims \ref{claim:class} and Corollary \ref{cor:boolean}, $X_{\phi}$ is a possibility domain if and only if it admits a binary non-dictatorial aggregator such that there exists a $d\in\{1,2\}$ and a non-empty $J\subseteq\{1,\ldots,k\}$, where, for all $j\in J\cup\{k+1,k+2,k+3\}$, $f_j=\pr_d^2$. It is not difficult to see that such an aggregator exists for $d=1$ if and only if it does for $d=2$. Thus, we can safely assume that $d=1$.

Before proving that our reduction works, we need some notation. For a domain $Y\subseteq\{0,1\}^m$ and a non-empty set of indices $I\subseteq\m$, let $Y_I$ be the projection of $Y$ to the indices of $I$, that is, the set of partial vectors with coordinates in $I$ that can be extended to vectors of $Y$. Also, for two domains $Y,Z\subseteq\{0,1\}^m$ of size $n$, viewed as $n\times m$ matrices, we write $Y\approx Z$ if, by permuting the columns of $Z$, we can obtain $Y$.

First, assume that there exist $1<l<k$ variables $x_{i_1},\ldots,x_{i_l}\in\{x_1,\ldots,x_k\}$ such that $\chi$ is independent from all of them. Let also $J=\{1,\ldots,k\}\setminus\{i_1,\ldots,i_l\}$. Then, \eqref{modx} can be written as:$$X_{\phi}\approx\{0,1\}^l\times\Bigg(\Big(\rMod(\chi)_J\times\{(0,0,0),(1,1,1)\}\Big)\cup\Big(\rMod(\neg\chi)_J\times\imp\Big)\Bigg).$$ It is straightforward to observe that any $m$-tuple $F=(f_1,\ldots,f_m)$ of binary functions, such that $f_{l+1}=\ldots=f_m=\pr_d^2$, $d=1,2$ is an aggregator for $X_{\phi}$. Thus $X_{\phi}$ is a possibility domain.

Now, assume that $X_{\phi}$ is a possibility domain and $F=(f_1,\ldots,f_m)$ a binary non-dictatorial aggregator for $X_{\phi}$. Let also $(I,J)$ be a partition of $\{1,\ldots,k\}$, such that, for all $j\in J\cup\{k+1,k+2,k+3\}$, $f_j=\pr_1^2$ and for all $i\in I$, $f_i\neq \pr_1^2$. Again, to simplify things, assume that there exists an $s\in\{1,\ldots,k-1\}$ such that $I=\{1,\ldots,s\}$, $J=\{s+1,\ldots,k\}$. We consider the following three cases:\begin{itemize}
    \item If $f_i:=\pr_2^2$, for all $i\in I$, we show that $\chi$ is independent from $x_1,\ldots,x_s$. Suppose there exist vectors $\bar{a},\bar{b}\in\{0,1\}^s$ and $\bar{c}\in\{0,1\}^{k-s}$, such that $(\bar{a},\bar{c})\in\rMod(\chi)$ and $(\bar{b},\bar{c})\in\rMod(\neg\chi)$. Then, $(\bar{a},\bar{c},0,0,0)\in X_{\phi}$ and $(\bar{a},\bar{c},0,0,1)\in X_{\phi}$. Since $F$ is an aggregator for $X_{\phi}$, it must hold that:$$F((\bar{a},\bar{c},0,0,0),(\bar{b},\bar{c},0,0,1))=(\bar{b},\bar{c},0,0,0)\in X_{\phi},$$ which, by \eqref{modx}, implies that $(\bar{b},\bar{c})\in\rMod(\chi)\cap\rMod(\neg\chi)$. Contradiction. Since $\bar{a},\bar{b}$ and $\bar{c}$ where chosen arbitrarily, it follows that $\chi$ is independent from $x_1,\ldots,x_s$.
     \item If $f_i$ is symmetric, for all $i\in I$, we show that $\chi$ is independent from $x_1,\ldots,x_s$. Suppose there exist vectors $\bar{a},\bar{b}\in\{0,1\}^s$ and $\bar{c}\in\{0,1\}^{k-s}$, such that $(\bar{a},\bar{c})\in\rMod(\chi)$ and $(\bar{b},\bar{c})\in\rMod(\neg\chi)$. Also, assume that $(f_1,\ldots,f_s)(\bar{a},\bar{b}):=\bar{z}$. Then, $(\bar{a},\bar{c},0,0,0)\in X_{\phi}$ and $(\bar{b},\bar{c},0,0,1)\in X_{\phi}$. Since $F$ is an aggregator for $X_{\phi}$, it must hold that:\begin{align*}F((\bar{a},\bar{c},0,0,0),(\bar{b},\bar{c},0,0,1)) & =(\bar{z},\bar{c},0,0,0)\in X_{\phi},\\ F((\bar{b},\bar{c},0,0,1),(\bar{a},\bar{c},0,0,0)) & =(\bar{z},\bar{c},0,0,1)\in X_{\phi},\end{align*} which, by \eqref{modx}, implies that $(\bar{z},\bar{c})\in\rMod(\chi)\cap\rMod(\neg\chi)$. Contradiction. Since $\bar{a}$, $\bar{b}$ and $\bar{c}$ where chosen arbitrarily, it follows that $\chi$ is independent from $x_1,\ldots,x_s$.
     \item If there is a partition $(I_1,I_2)$ of $I$, such that $f_i=\pr_2^2$ for all $i\in I_1$ and $f_i$ is symmetric for all $i\in I_2$, we show that $\chi$ is independent from all $x_i$ such that $i\in I_1$. Assume again w.l.o.g. that there is a $t\in\{1,\ldots,s-1\}$ such that $I_1=\{1,\ldots,t\}$ and $I_2=\{t+1,\ldots,s\}$. Suppose there exist vectors $\bar{a},\bar{b}\in\{0,1\}^t$ and $\bar{c}\in\{0,1\}^{k-t}$, such that $(\bar{a},\bar{c})\in\rMod(\chi)$ and $(\bar{b},\bar{c})\in\rMod(\neg\chi)$. Then, $(\bar{a},\bar{c},0,0,0)\in X_{\phi}$ and $(\bar{a},\bar{c},0,0,1)\in X_{\phi}$. Since $F$ is an aggregator for $X_{\phi}$, it must hold that:$$F((\bar{a},\bar{c},0,0,0),(\bar{b},\bar{c},0,0,1))=(\bar{b},\bar{c},0,0,0)\in X_{\phi},$$ which, by \eqref{modx}, implies that $(\bar{b},\bar{c})\in\rMod(\chi)\cap\rMod(\neg\chi)$. Contradiction. Since $\bar{a},\bar{b}$ and $\bar{c}$ where chosen arbitrarily, it follows that $\phi$ is independent from $x_1,\ldots,x_t$.
\end{itemize}
This concludes the proof of both the reduction and Theorem \ref{thm:possintcon}.\hfill$\Box$\vspace{0.2cm}

We now show the analogous result for local possibility domains. Fortunately, the proofs of the bounds here are relatively shorter.

\begin{theorem}\label{thm:locposintcon}
Deciding, on input $\phi$,  whether or not $X_{\phi}$  admits a locally non-dictatorial aggregator is (i) in $\rS^P_2\cap\rP^P_2$ and (ii) $\coNP\text{-hard}$.
\end{theorem}
\textbf{Proof}
(i) We follow the proof of Theorem \ref{thm:possintcon}. By Corollary \ref{cor:lpdchar} we have only three functions, namely $\wedge^{(3)},\vee^{(3)},\oplus$, which, when combined to an $m$-ary tuple $F=(f_1,\ldots,f_m)$, form a local non-dictatorial aggregator for $X$. Thus, the proof is exactly the same, with the difference that we now have $3^m$ tuples that can be encoded as $(6\times3)m$-binary sequences and we conclude that the problem is in $\rS^P_2$.

For the containment in $\rP^P_2$, we argue as follows. Given access to $H_{X_{\phi}}$ we can, by Corollary \ref{cor:tractlpd}, obtain in polynomial time a ternary WNU aggregator $G=(g_1,\ldots,g_m)$ such that $X_{\phi}$ is a local possibility domain if and only if it admits $G$. Whether $X_{\phi}$ admits $G$ or not is in $\coNP$ (in the same way we check if $X_{\phi}$ is affine) and thus in $\rP_2^P$ too. Since the size of $H_{X_{\phi}}$ is polynomial in that of $\phi$, it suffices to show that testing whether there is no edge $uu'_i\rightarrow\rightarrow vv'_j$ in $H_{X_{\phi}}$ is in $\rP_2^P$. This problem can be expressed as: \begin{quote}
For all assignments $\ba=(a_1,\ldots,a_m)$, $\bb=(b_1,\ldots,b_m)$, there exists an assignment $\bc=(c_1,\ldots,c_m)$ such that, if $\ba$, $\bb$ satisfy $\phi$ and $a_i=u$, $b_i=u'$, $a_j=v$, $b_j=v'$, then $\bc$ satisfies $\phi$, $c_i=u$, $c_j=v'$ and $c_k\in\{a_k,b_k\}$ for all $k\in\m\setminus\{i,j\}$.
\end{quote} Thus, the proof is complete.

(ii) We show that the problem of whether a logical formula $\chi$, defined on $k$ variables $x_1,\ldots,x_k$, is unsatisfiable, reduces to that of deciding if the truth set of a formula is a local possibility domain.

Let $\psi(y_1,y_2,y_3)$ be the propositional formula with $\rMod(\psi)=\imp$, where $\{y_1,y_2,y_3\}\cap\{x_1,\ldots,x_k\}=\emptyset$. Consider the formula: $$\phi=(\chi(x_1,\ldots.x_k)\wedge\psi(y_1,y_2,y_3))\vee(z\rightarrow w),$$ where $z$ and $w$ are variables not among those of $\chi$ or $\psi$. First note that the length of $\phi$ is again linear to that of $\chi$, since $|\phi|=|\chi|+8$ and thus the construction is polynomial.

Let $m:=k+5$. We prove the following claim.
\begin{claim}\label{claim:unsat}
$\chi$ is unsatisfiable if and only if \begin{equation}\label{modx2} X_{\phi}=\Big(\rMod(\chi)\times\imp\times\{0,1\}^2\Big)\cup\Big(\{0,1\}^{k+3}\times\{(0,0),(0,1),(1,1)\}\Big),\end{equation} is a local possibility domain.\end{claim}
\textit{Proof of Claim:}
First, assume $\chi$ is unsatisfiable. Then, \eqref{modx2}:$$X_{\phi}=\{0,1\}^{k+3}\times\{(0,0),(0,1),(1,1)\},$$ which is a local possibility domain, since it admits for example the binary aggregator $F=(f_1,\ldots,f_m)$, where $f_j=\wedge$, $j=1,\ldots,m$.

On the other hand, let $\rMod(\phi)$ be a local possibility domain, and assume $\chi$ is satisfiable by some assignment $\bar{a}=(a_1,\ldots,a_k)$. Since $\rMod(\phi)$ is a local possibility domain, by Theorem \ref{thm:updcar}, it admits a ternary locally non-dictatorial aggregator $F=(f_1,\ldots,f_m)$.

By Lemma \ref{lem:imp}, $\imp$ is an impossibility domain. Thus, there exist $\bar{b}^i=(b_1^i,b_2^i,b_3^i)\in\imp$, $i=1,2,3$, such that:$$\bar{c}:=(f_{k+1},f_{k+2},f_{k+3})(\bar{b}^1,\bar{b}^2,\bar{b}^3)\notin\imp.$$ Since $\bar{b}^1,\bar{b}^2,\bar{b}^3\in\imp$, it holds that $(\bar{a},\bar{b}^i,1,0)\in X_{\phi}$, for $i=1,2,3$. On the other hand:$$F((\bar{a},\bar{b}^1,1,0),(\bar{a},\bar{b}^2,1,0),(\bar{a},\bar{b}^3,1,0))=(\bar{a},\bar{c},1,0)\notin X_{\phi}.$$ Thus $F$ is not an aggregator for $X_{\phi}$. Contradiction.\hfill$\Box$

This concludes the proof of Theorem \ref{thm:locposintcon}.

\subsection{Agendas}\label{ssec:agenda}
Suppose now that we have an agenda $\ph=(\phi_1,\ldots,\phi_m)$. We prove the following upper and lower bounds to the complexity of deciding whether $X_{\ph}$ is a (local) possibility domain. For upper bounds, we use a version of Lemma \ref{lem:closure} tailored for agendas.
\begin{lemma}\label{lem:agclosure}
Let $F=(f_1,\ldots,f_m)$ be an $m$-tuple of $n$-ary conservative and polynomial-time computable functions. Deciding, given the agenda $\ph=(\phi_1,\ldots,\phi_m)$, whether $F$ is an aggregator for $X_{\ph}$ is in $\rP_2^P$.	
\end{lemma}
\textbf{Proof} The result is immediate since the problem can be cast as follows:
   \begin{quote}
    For all $m$-tuples $\x^1,\ldots,\x^n\in\{0,1\}^m$, if $\bigwedge_{j=1}^m\phi_j^{x^1_j},\ldots,\bigwedge_{j=1}^m\phi_j^{x^n_j}$ are all satisfiable, then so is $\bigwedge_{j=1}^m\phi_j^{w_j}$, where $w_j=F(x^1_j,\ldots,x^n_j)$, $j=1,\ldots,m$
\end{quote}
and since checking if a formula is satisfied by a specific assignment can be done in polynomial time.\hfill$\Box$\vspace{0.2cm}

In terms of lower bounds, a straightforward idea would be to construct, given an integrity constraint, an agenda with the same domain, since that would immediately imply that the lower bounds for the integrity constraints carry on to the agendas. Unfortunately, as discussed above, this is an $\FNP$-complete problem. However, by Endriss et al. \cite[Proposition 3]{endriss2016succinctness}, we can obtain the following result.
\begin{corollary}\label{cor:succ}
Given an integrity constraint $\phi$ and a satisfying assignment $\bar{a}$ of $\phi$, we can construct an agenda $\bar{\phi}$, of polynomial size in the length of $\phi$, such that $X_{\ph}=X_{\phi}$. Furthermore, this construction is polynomial in the size of $\phi$.
\end{corollary}
Corollary \ref{cor:succ} follows directly from the proof of \cite[Proposition 3]{endriss2016succinctness}. The reason this result does not imply the existence of a polynomial reduction, is that to construct $\ph$, one needs a satisfying assignment of $\phi$. And of course, finding such an assignment is intractable. Fortunately, we can get past that in the problems we consider.

\begin{theorem}\label{pdagenda}
Given the agenda $\ph=(\phi_1,\ldots,\phi_m)$, the question whether $X_{\ph}$ admits a non-dictatorial aggregator is (i) in $\rD^P_3$ and (ii) $\coNP$-hard.
\end{theorem}
\textbf{Proof}
(i) We will show that the problem can be decided in $P$ with an oracle in $\rS_2^P$. By Corollary \ref{cor:boolean}, $X_{\ph}$ is a possibility domain if and only if it is affine or it admits a binary non-dictatorial aggregator. The problem of whether $X_{\ph}$ is affine is in $\rP_2^P\subseteq\rD_3^P$ by Lemma \ref{lem:agclosure}, using $\bar{\oplus}$ as the aggregator.

It remains to show that the latter problem, or equivalently the problem of checking if $H_{X_{\ph}}$ is strongly connected, is in $\rD_3^P$. Since $H_{X_{\ph}}$ has $2m$ vertices, its size is polynomial in that of the agenda. Also, checking if a graph is strongly connected is in $\P$. Thus it suffices to show that $H_{X_{\ph}}$ can be constructed within polynomial time with an oracle in $\rS_2^P$.

\begin{claim}\label{HXconstr}
Given the agenda $\ph=(\phi_1,\ldots,\phi_m)$, the graph $H_{X_{\ph}}$ can be constructed within polynomial time in the size of $\ph$ with an oracle in $\rS_2^P$.
\end{claim}
\textit{Proof of Claim:}
Consider two vertices $uu'_s$ and $vv'_t$ of $H_{X_{\ph}}$, with $s\neq t$. To decide if there is an edge from $uu'_s$ to $vv'_t$, is suffices to check the following:
\begin{quote}
    There exist binary $m$-sequences $\x,\y$ such that:\begin{itemize}
        \item both $\bigwedge_{j=1}^m\phi_j^{x_j}$ and $\bigwedge_{j=1}^m\phi_j^{y_j}$ are satisfiable,
        \item $x_s=u$, $x_t=v$, $y_s=u'$ and $y_t=v'$ and
        \item for all $m$-sequences $\bar{z}$, either $\bigwedge_{j=1}^m\phi_j^{z_j}$ is not satisfiable or at least one of the following is not true: (i) $z_k=u$, (ii) $z_l=v'$, (iii) $z_j\in\{x_j,y_j\}$ for $j=1,\ldots,m$.
    \end{itemize}
\end{quote}
Notice that this can be done with an oracle in $\rS_2^P$.\hfill$\Box$

(ii) Let $\chi(x_1,\ldots,x_k)$ be a propositional formula on $k$ variables. We construct, in polynomial time, an agenda $\ph$ such that $\chi$ is independent from at least one of its variables if and only if $X_{\ph}$ is a possibility domain. In all that follows, $m=k+3$.

Let: $$\phi(x_1,\ldots,x_k,y_1,y_2,y_3)=\chi(x_1,\ldots,x_k)\oplus\psi(y_1,y_2,y_3),$$ where: $\{x_1,\ldots,x_l\}$ and $\{y_1,y_2,y_3\}$ are disjoint sets of variables and $\rMod(\psi)=\imp$. As in Theorem \ref{thm:possintcon}, we have that $\phi$'s length is linear to that of $\chi$ and that:
\begin{equation}\label{modxx}
X_{\phi}=\Big(\rMod(\chi)\times\{(0,0,0),(1,1,1)\}\Big)\cup\Big(\rMod(\neg\chi)\times\imp\Big).\end{equation}
Pick an arbitrary vector $\bar{a}=(a_1,\ldots,a_k)\in\{0,1\}^k$ and set:\begin{equation*}\bar{b}=(b_1,\ldots,b_m):=\begin{cases} (a_1,\ldots,a_k,0,0,0) & \text{if }\chi(a_1,\ldots,a_k)=1,\\
(a_1,\ldots,a_k,0,0,1) & \text{ else.}\end{cases}
\end{equation*}
In both cases, $\bar{b}$ satisfies $\phi$. Thus, we can use Corollary \ref{cor:succ}, to construct an agenda $\ph$ such that $X_{\ph}=X_{\phi}$, whose size is polynomial in the length of $\phi$ and thus in that of $\chi$ too. Since deciding whether $\bar{a}$ satisfies $\chi$ or not can be done in polynomial time, our construction is polynomial. Also, in Claim \ref{claim:indvar} of Theorem \ref{thm:possintcon}, we prove that $\chi$ is independent from at least one of its variables if and only if $X_{\phi}$, and thus $X_{\ph}$ too, is a possibility domain. \hfill $\Box$\vspace{0.2cm}

In what concerns local possibility domains, we have the analogous bounds.
\begin{theorem}\label{agendalocal}
Given the agenda $\ph=(\phi_1,\ldots,\phi_m)$, the question of whether $X_{\ph}$ admits a locally non-dictatorial aggregator is (i) in $\rD^P_3$ and (ii) $\coNP$-hard.
\end{theorem}
\textbf{Proof}
(i) We have already argued in Claim \ref{HXconstr} that $H_{X_{\bar{\phi}}}$ can be constructed in polynomial time with an oracle in $\rS_2^P$. Then, by Corollary \ref{HXmaxsym} we can obtain, in polynomial time, a maximum symmetric aggregator $F=(f_1,\ldots,f_m)$ for $X$. We can easily now obtain the aggregator $G=(g_1,\ldots,g_m)$ for $X$ of the first paragraph Corollary \ref{cor:tractlpd}, where $g_j\in\{\wedge^{(3)},\vee^{(3)},\oplus\}$, $j=1,\ldots,m$.
%phk and
Furthermore, testing whether $X_{\ph}$ admits this aggregator  is in $\rP_2^P$ by Lemma \ref{lem:agclosure}.

(ii) Let $\chi(x_1,\ldots,x_k)$ be a propositional formula. We construct an agenda $\ph$ such that $X_{\ph}$ is a local possibility domain if and only if $\chi$ is unsatisfiable. Let $\psi(y_1,y_2,y_3)$ be again the propositional formula with $\rMod(\psi)=\imp$, where $\{y_1,y_2,y_3\}\cap\{x_1,\ldots,x_k\}=\emptyset$. Consider the formula: $$\phi=(\chi(x_1,\ldots.x_k)\wedge\psi(y_1,y_2,y_3))\vee(z\rightarrow w),$$ where $z$ and $w$ are variables not among those of $\chi$ or $\psi$.

Let $m:=k+5$. By \eqref{modx2} we have that: \begin{equation}\label{modxx2} X_{\phi}=\Big(\rMod(\chi)\times\imp\times\{0,1\}^2\Big)\cup\Big(\{0,1\}^{k+3}\times\{(0,0),(0,1),(1,1)\}\Big).\end{equation}

Again, pick an arbitrary vector $\bar{a}=(a_1,\ldots,a_k)\in\{0,1\}^k$ and set:\begin{equation*}\bar{b}=(b_1,\ldots,b_m):=(a_1,\ldots,a_k,0,0,0,0,0).
\end{equation*}
$\bar{b}$ satisfies $\phi$. Thus, we can use Proposition $3$ of \cite{endriss2016succinctness}, to construct an agenda $\ph$ such that $X_{\ph}=X_{\phi}$, whose size is polynomial in the length of $\phi$ and thus in that of $\chi$ too. Also, in Claim \ref{claim:unsat}, we prove that $\chi$ is unsatisfiable if and only if $X_{\phi}$, and thus $X_{\ph}$ too, is a local possibility domain. \hfill $\Box$.\vspace{0.2cm}

A related result, that has been answered in \cite{endriss2012complexity}, is whether the domain of an agenda admits the majority aggregator. In \cite{nehring2007structure} and \cite{endriss2012complexity}, such agendas are characterized as those satisfying the \emph{median property}, that is, agendas whose every inconsistent subset, contains an inconsistent subset of size $2$. Endriss et al. \cite{endriss2012complexity} show that checking if an agenda satisfies the median property is $\rP_2^P$-complete. Unfortunately, this does not extend to the problem of determining if the domain of an agenda is a possibility domain, since, even though checking for the minority aggregator is in $\rP_2^P$, the existence of binary non-dictatorial aggregators seems to be a harder problem.

\begin{remark}
It is clear that the upper and lower complexity bounds obtained here  do not match. By Theorem \ref{pdagenda} and Theorem  \ref{agendalocal}, it is also clear that the gap between the upper and lower bounds in the case of the agendas comes from the fact that the lower bounds are obtained via the integrity constraints.

What is not so clear, however, is why in the case of integrity constraints the upper bound is $\rS_2^P\cap\rP_2^P$ and in the case of the agendas the upper bound is $\rD_3^P$, whereas one might except a full computational jump in the \PH, or for the bounds to be the same. First, note that,
in the case of the integrity constraints, it does not seem plausible  to construct $H_{X_{\phi}}$ in polynomial time, with access only to an \NP \ oracle. This is so because, to decide if there is an edge $uu'_s\rightarrow vv'_t$ in $H_X$, one needs to ``guess'' two satisfying assignments for $\phi$, such that there is no other satisfying assignment with certain properties. This naturally points to an existential and a universal certificate, something that an \NP \ oracle alone cannot provide. On the other hand, for the same reason, it seems improbable that finding an aggregator $F$ for $X_{\ph}$ can be in $\Sigma_2^P$, since we need to guess a candidate $F$, such that for all vectors of $0/1$ judgments for the formulas of the agenda, we can decide the satisfiability of a specific formula. This again needs more certificates than $\Sigma_2^P$ (or $\Pi_2^P$) can provide.

The above arguments justify only why our proofs do not yield any better complexity bounds. Naturally, it is conceivable that other approaches may yield tighter bounds.\hfill$\diamond$
\end{remark}

\subsection{Other Types of Non-dictatorial Aggregation}\label{ssec:other}
In this subsection, we quickly extend the results of the previous subsections in four cases of non-dictatorial aggregation that have been used in the literature. Namely, we discuss \emph{generalized dictatorships}, \emph{anonymous}, \emph{monotone} and \emph{systematic} aggregators. We only consider the case where we search if an implicitly given Boolean domain admits such aggregators. The case where the domain is given explicitly has been shown to be tractable by Diaz et al. \cite{diaz2019syntactic}, since the results there directly provide the required aggregators.

\subsubsection{Generalized Dictatorships}
We begin by defining generalized dictatorships, that is, aggregators whose output is always a vector of the input.
\begin{definition}\label{def:gendict}
Let $F =(f_1,\ldots,f_m)$ be an $m$-tuple of $n$-ary conservative functions. $F$ is a \emph{generalized dictatorship for a domain} $X\subseteq\{0,1\}^m$, if, for any $x^1,\ldots,x^n\in X$, it holds that:
\begin{equation}\label{gendic}F(x^1,\ldots,x^n):=(f_1(x_1),\ldots,f_m(x_m))\in\{x^1,\ldots,x^n\}.\end{equation}
\end{definition}
Obviously, generalized dictatorships of a domain $X$ are always aggregators for $X$. Generalized dictatorships were first defined by Grandi and Endriss \cite{grandi2013lifting}, where they were defined independently of a specific domain. That is, Eq. \eqref{gendic} holds for all $x^1,\ldots,x^n\in\{0,1\}^m$ in their context. With this stronger definition, Grandi and Endriss \cite[Theorem 16]{grandi2013lifting} show that the class of generalized dictators coincides with that of functions that are aggregators for every domain $X\subseteq\{0,1\}^m$.

Generalized dictatorships select a possibly different dictator for each input. All dictatorial aggregators are obviously generalized dictatorships. Depending on the domain, different kinds of non-dictatorial aggregators can be generalized dictatorships. For example, $\overline{\maj}$ and $\bar{\oplus}$ are both generalized dictatorships for any domain with only two elements. In fact, any domain with only two elements admits only generalized dictatorships as aggregators. The following characterization has been proven in \cite{diaz2019syntactic}.

\begin{customthm}{F}{\cite[Theorem 11]{diaz2019syntactic}}\label{thm:gendict}
A domain $X\subseteq\{0,1\}^m$ admits an aggregator that is not a generalized dictatorship if and only if it is a possibility domain with at least three elements.
\end{customthm}

Using that, we can easily prove the same complexity bounds we had for deciding if the domain of an integrity constraint admits a non-dictatorial aggregator.

\begin{corollary}\label{cor:gendictup}
Deciding, on input $\phi$,  whether $X_{\phi}$  admits an aggregator that is not a generalized dictatorship is (i) in $\rS^P_2\cap\rP_2^P$ and (ii) $\coNP$-hard.	
\end{corollary}
\textbf{Proof} By Theorem \ref{thm:gendict}, it suffices to show that deciding if $X_{\phi}$ has at least three elements is in $\rS^P_2\cap\rP_2^P$. Indeed, this can be written as:
\begin{quote}
    There exist tuples $\x,\y,\bar{z}\in\{0,1\}^m$ that are pairwise distinct, such that all three satisfy $\phi$.
\end{quote}
Thus, deciding if $|X_{\phi}|\geq 3$ is in $\rS_1^P\subseteq\rS^P_2\cap\rP_2^P$. The rest of the proof is identical with that of Theorem \ref{thm:possintcon}.

(ii) Immediate by Theorems \ref{thm:gendict} and \ref{thm:possintcon}, since the domain of \eqref{modx} has more than two elements.\hfill$\Box$\vspace{0.2cm}

In case the domain is provided via an agenda, we can again easily obtain the same bounds.
\begin{corollary}\label{cor:gendictagendaup}
Given the agenda $\ph=(\phi_1,\ldots,\phi_m)$, the question whether $X_{\ph}$ admits an aggregator that is not a generalized dictatorship is (i) in $\rD^P_3$ and (ii) $\coNP$-hard.
\end{corollary}
\textbf{Proof}
(i) Again by Theorem \ref{thm:gendict}, it suffices to show that deciding whether $X_{\ph}$ has at least three elements is in $\rD_3^P$. Indeed, the problem can be written as:\begin{quote}
    There exist $m$-tuples $\x,\y,\bar{z}\in\{0,1\}^m$ which are pairwise distinct and such that $\bigwedge_{j=1}^m\phi_j^{x_j}$, $\bigwedge_{j=1}^m\phi_j^{y_j}$ and $\bigwedge_{j=1}^m\phi_j^{z_j}$ are all satisfiable.
\end{quote}
Thus, it is in $\rS_2^P\subseteq\rD_3^P$. The rest of the proof is identical that of Theorem \ref{pdagenda}.

(ii) Immediate by Theorems \ref{thm:gendict} and \ref{pdagenda}, since the domain in \eqref{modxx} has more than two elements.\hfill$\Box$\vspace{0.2cm}

\subsubsection{Anonymous aggregators}
An aggregator is anonymous if it is not affected by permutations of its input. Such an aggregator treats each member of the population ``fairly'', in the sense that the outcome of the aggregation depends only on the votes that have been cast and not on which person selected which options.
\begin{definition}\label{def:anonymous}
Let $X\subseteq\{0,1\}^m$. An $n$-ary aggregator $F=(f_1,\ldots,f_m)$ for $X$ is anonymous, if it holds that for all $j\in\m$ and for any permutation $p:\{1,\ldots,n\}\mapsto\{1,\ldots,n\}$:$$f_j(a_1,\ldots,a_n)=f_j(a_{p(1)},\ldots,a_{p(n)}),$$ for all $a_1,\ldots,a_n\in\{0,1\}$.
\end{definition}
It is immediate to observe that a ternary such aggregator is always WNU. Thus, using Theorem \ref{thm:updcar}, Kirousis et al. proved the following result.
\begin{corollary}{\cite[Corollary 5.11]{kirousis2019aggregation}}\label{cor:anon} $X\subseteq\{0,1\}^m$ is a local possibility domain if and only if it admits an anonymous aggregator.
\end{corollary}
In fact, a local possibility domain always admits a ternary such aggregator. We can now obtain all the complexity bounds we have for local non-dictatorial aggregators, in the case where we search for anonymous ones.
\begin{corollary}\label{results:anon}
\begin{itemize}
    \item[A.] Deciding, on input $\phi$,  whether $X_{\phi}$  admits an anonymous aggregator is i. in $\rS_2^P\cap\rP_2^P$ and ii. $\coNP$-hard.	
    \item[B.] Given the agenda $\ph=(\phi_1,\ldots,\phi_m)$, the question whether $X_{\ph}$ admits an anonymous aggregator is i. in $\rD_3^P$ and ii. $\coNP$-hard.
\end{itemize}
\end{corollary}
\textbf{Proof} \begin{itemize}
    \item[Ai.] Immediate by Corollary \ref{cor:anon} and Theorem \ref{thm:locposintcon}.
    \item[Aii.] Immediate by Corollary \ref{cor:anon} and Theorem \ref{thm:locposintcon}.
    \item[Bi.] Immediate by Corollary \ref{cor:anon} and Theorem \ref{agendalocal}.
    \item[Bii.] Immediate by Corollary \ref{cor:anon} and Theorem \ref{agendalocal}.\hfill$\Box$\end{itemize}

\subsubsection{Monotone aggregators}
Monotone aggregators are aggregators whose output does not change if more individuals agree with it. This is one of the features that makes majority voting appealing.
\begin{definition}\label{def:monotone}
Let $X\subseteq\{0,1\}^m$. An $n$-ary aggregator $F=(f_1,\ldots,f_m)$ for $X$ is monotone, if it holds that for all $j\in\m$ and for all $i\in\{1,\ldots,n\}$: $$f_j(a_1,\ldots,a_{i-1},0,a_{i+1},\ldots,a_n)=1\Rightarrow f_j(a_1,\ldots,a_{i-1},1,a_{i+1},\ldots,a_n)=1.$$
\end{definition}
It is not difficult to see that all binary aggregators (both dictatorial and non-dictatorial) have that property. Also, in \cite{dokow2009aggregation} and \cite{kirousis2019aggregation}, it has been proven that if a domain admits a majority aggregator, it also admits a binary non-dictatorial one. Combining this with Corollary \ref{cor:boolean}, we obtain the following result.
\begin{customthm}{G}{\cite[Theorem $12$]{diaz2019syntactic}}\label{thm:monotone}
A domain $X\subseteq\{0,1\}^m$ admits a monotone non-dictatorial aggregator of some arity if and only if it admits a binary non-dictatorial one.
\end{customthm}
Thus we can again easily obtain the required complexity bounds.
\begin{corollary}\label{results:monotone}
\begin{itemize}
    \item[A.] Deciding, on input $\phi$,  whether $X_{\phi}$  admits a monotone non-dictatorial aggregator is i. in $\rS_2^P\cap\rP_2^P$ and ii. $\coNP$-hard.
    \item[B.] Given the agenda $\ph=(\phi_1,\ldots,\phi_m)$, the question whether $X_{\ph}$ admits a monotone non-dictatorial aggregator is i. in $\rD_3^P$ and ii. $\coNP$-hard.
\end{itemize}
\end{corollary}
\textbf{Proof} \begin{itemize}
    \item[Ai.] Immediate by Theorem \ref{thm:monotone} and Theorem \ref{thm:possintcon}.
    \item[Aii.] Immediate by Theorem \ref{thm:monotone} and by noticing that the only available aggregators for the construction of \eqref{modx} in Theorem \ref{thm:locposintcon} are binary non-dictatorial ones.
    \item[Bi.] Immediate by Theorem \ref{thm:monotone} and Theorem \ref{pdagenda}.
    \item[Bii.] Immediate by Theorem \ref{thm:monotone} and by noticing that the only available aggregators for the construction of \eqref{modxx} in Theorem \ref{pdagenda} are binary non-dictatorial ones.\hfill$\Box$\end{itemize}

\subsubsection{Systematic aggregators}
Systematicity is the requirement that an aggregator treats all issues in the same way. This can be a very natural requirement for a set of issues of the same nature. We have already seen that systematic aggregators of a domain correspond to polymorphisms. Recall that we denote a systematic aggregator $F=(f_1,\ldots,f_m)$, where $f_1=\cdots=f_m=f$, by $\f$.

Polymorphisms have been extensively studied in the literature and they play a central role in Post's \cite{post1941two} seminal work, where he provides a full classification of Boolean clones, that is sets of Boolean functions closed under composition. Post's work can provide results in Complexity Theory too; see for example \cite{bohler2003playing,bohler2004playing,jeavons1995algebraic,jeavons1999determine}. Here, we use a result concerning the unanimous polymorphisms (equivalently, the systematic aggregators), that a Boolean domain admits. We say that an $n$-ary function $f$ is \emph{essentially unary}, if there exists a \emph{unary} Boolean function $g$ and an $i\in\{1,\ldots,n\}$, such that: $$f(x_1,\ldots,x_n)=g(x_i),$$ for all $x_1,\ldots,x_n\in\{0,1\}$. Obviously, the only unanimous such functions are the projections.
\begin{corollary}\label{cor:pol}
Let $X\subseteq\{0,1\}^n$ be a Boolean domain. Then, either $X$ admits only essentially unary functions, or it is closed under $\wedge$, $\vee$, $\maj$ or $\oplus$.
\end{corollary}
This result can be obtained directly by Post's Lattice, without considering complexity theoretic notions. For a direct algebraic approach, see also \cite[Proposition $1.12$]{szendrei1986clones} (by noting that the only Boolean \emph{semi-projections} of arity at least $3$ are projections).

Corollary \ref{cor:pol} translates in our framework as follows.
\begin{corollary}\label{cor:system}
Let $X\subseteq\{0,1\}^n$ be a Boolean domain. Then $X$ admits a systematic non-dictatorial aggregator if and only if it admits the aggregators $\bar{\wedge}$, $\bar{\vee}$, $\overline{\maj}$ or $\bar{\oplus}$).
\end{corollary}

In case of integrity constraints, the problem of detecting if their domains admit systematic non-dictatorial aggregators is $\coNP$-complete.

\begin{proposition}\label{icshaefer}
Deciding, on input $\phi$,  whether $X_{\phi}$  admits a systematic non-dictatorial aggregator is $\coNP$-complete.
\end{proposition}
\textbf{Proof} Checking closure under $\wedge$, $\vee$, $\maj$ and $\oplus$ is in $\coNP$ by Lemma \ref{lem:closure}. Thus we only need to show $\coNP$-hardness. We reduce from the known $\coNP$-complete problem of \emph{tautology}, where we check if a propositional formula is satisfied by all assignments of values.

Let $\chi(x_1,\ldots,x_k)$ be the input propositional formula on $k$ variables, and $\psi(y_1,y_2,y_3)$ be the formula such that $\rMod(\psi)=\imp$, where $\{x_1,\ldots,x_m\}\cap\{y_1,y_2,y_3\}=\emptyset$. Let also $m=k+3$.

Consider the formula: $$\phi(x_1,\ldots,x_m,y_1,y_2,y_3)=\chi\vee\psi.$$ If $\chi$ is a tautology, $X_{\phi}=\{0,1\}^m$, which is closed under $\wedge$, $\vee$, $\maj$ and $\oplus$. Otherwise: $$X_{\phi}=\Big(\rMod(\chi)\times\{0,1\}^3\Big)\cup\Big(\rMod(\neg\chi)\times\imp\Big).$$ Easily now, if $\bar{a}=(a_1,\ldots,a_k)$ does not satisfy $\phi$, and $f\in\{\wedge,\vee,\maj,\oplus\}$, it holds that $(\bar{a},0,0,1)$, $(\bar{a},0,1,0)$, $(\bar{a},1,0,0)\in X_{\phi}$, but:$$\bar{b}:=\f((\bar{a},1,0,0),(\bar{a},0,1,0),(\bar{a},1,0,0))=\{(\bar{a},0,0,0),(\bar{a},1,1,1)\}.$$ Thus, since $\bar{b}\notin X_{\phi}$, $X_{\phi}$ does not admit any systematic non-dictatorial aggregator.\hfill$\Box$\vspace{0.2cm}
Finally, we can obtain the corresponding results in the case the domain is provided via an agenda.
\begin{proposition}\label{agendashaefer}
Deciding, on input $\ph=(\phi_1,\ldots,\phi_m)$,  whether $X_{\ph}$  admits a systematic non-dictatorial aggregator is in $\rP_2^P$ and $\coNP$-hard.
\end{proposition}
\textbf{Proof} By Lemma \ref{lem:agclosure}, we have membership in $\rP_2^P$ for detecting closure under $\wedge$, $\vee$, $\maj$ and $\oplus$. For $\coNP$-hardness, observe that given a formula $\chi$, we can again set $\phi=\chi\vee\psi$ and construct an agenda whose domain is the same with $X_{\phi}$ of Proposition \ref{icshaefer} in polynomial time, by Endriss et al. \cite[Proposition $3$]{endriss2016succinctness}. Thus, the same reduction as in Proposition \ref{icshaefer} works.\hfill $\Box$\vspace{0.2cm}

\section{Concluding Remarks}
In this paper, we established the first results concerning the tractability of non-dictatorial aggregation.
Specifically, we gave polynomial-time algorithms that take as input a  set $X$ of feasible evaluations and determine whether or not  $X$ is a possibility domain and a uniform possibility domain, respectively. In these algorithms,  the set $X$ of feasible evaluations is given to us explicitly, i.e., $X$ is given by listing all its elements. We also provided upper and lower bounds for the complexity of deciding if a domain $X$ is a possibility or a local possibility domain, in case it is given implicitly, either via an integrity constraint or by an agenda. Finally, we extended these results to other types of commonly used non-dictatorial aggregators. These bounds are preliminary, in the sense that they are not tight, especially in case the domains are provided via agendas. This merits further investigation, because implicitly given domains occur frequently in the field of judgment aggregation.

The work reported here assumes that the aggregators are conservative, an assumption that has been used heavily throughout the paper. There is a related, but weaker, notion of an \emph{idempotent} (or \emph{Paretian}) aggregator $F=(f_1,\ldots,f_m)$ where each $f_j$ is assumed to be an idempotent function, i.e., for all $x\in X_j$, we have that $f(x,\ldots,x) = x$. Clearly, every conservative aggregator is idempotent. In the Boolean framework, idempotent aggregators are conservative, but, in the non-Boolean framework, this  need not hold. It remains an open problem to investigate the computational complexity of the existence of non-dictatorial idempotent aggregators in the non-Boolean framework.

Another line of research could be to consider the aggregators in Subsection \ref{ssec:other} in the non-Boolean framework. This changes things in a non-trivial way. For example, Theorem \ref{thm:gendict} does not hold in the non-Boolean case, since for example a ternary minority aggregator can be a generalized dictatorship for specific domains. Finally, it would be interesting to follow along the lines of Terzopoulou et al. \cite{terzopoulou2018aggregating}, and consider settings where the individuals are allowed to cast votes on subsets of issues, both in the abstract and in the integrity constraint framework.\vspace{0.2cm}

\paragraph{{\bf Acknowledgments}} The research of Lefteris Kirousis was partially supported by the Special Account for Research Grants of the National and Kapodistrian University of Athens. The work of Phokion G.\ Kolaitis is partially supported by NSF Grant IIS-1814152.

%\bibliographystyle{plain}
%\bibliography{Compl_aggregation}

\begin{thebibliography}{10}

\bibitem{arrow1951social}
Kenneth~J. Arrow.
\newblock {\em Social choice and individual values}.
\newblock Wiley, New York, 1951.

\bibitem{bessiere2013detecting}
Christian Bessiere, Cl{\'e}ment Carbonnel, Emmanuel Hebrard, George Katsirelos,
  and Toby Walsh.
\newblock Detecting and exploiting subproblem tractability.
\newblock In {\em IJCAI}, pages 468--474, 2013.

\bibitem{bohler2003playing}
Elmar B{\"o}hler, Nadia Creignou, Steffen Reith, and Heribert Vollmer.
\newblock Playing with {B}oolean blocks, part {I}: Post’s lattice with
  applications to complexity theory.
\newblock In {\em ACM SIGACT Newsletter}, 2003.

\bibitem{bohler2004playing}
Elmar B{\"o}hler, Nadia Creignou, Steffen Reith, and Heribert Vollmer.
\newblock Playing with {B}oolean blocks, part {II}: Constraint satisfaction
  problems.
\newblock In {\em ACM SIGACT Newsletter}, 2004.

\bibitem{bulatov2006dichotomy}
Andrei~A. Bulatov.
\newblock A dichotomy theorem for constraint satisfaction problems on a
  3-element set.
\newblock {\em Journal of the ACM (JACM)}, 53(1):66--120, 2006.

\bibitem{bulatov2011complexity}
Andrei~A. Bulatov.
\newblock Complexity of conservative constraint satisfaction problems.
\newblock {\em ACM Transactions on Computational Logic (TOCL)}, 12(4):24, 2011.

\bibitem{carbonnel2016dichotomy}
Cl{\'e}ment Carbonnel.
\newblock The dichotomy for conservative constraint satisfaction is
  polynomially decidable.
\newblock In {\em International Conference on Principles and Practice of
  Constraint Programming}, pages 130--146. Springer, 2016.

\bibitem{carbonnel2016meta}
Cl{\'e}ment Carbonnel.
\newblock The meta-problem for conservative {M}al'tsev constraints.
\newblock In {\em Thirtieth AAAI Conference on Artificial Intelligence
  (AAAI-16)}, 2016.

\bibitem{diaz2019syntactic}
Josep Di\'az, Lefteris Kirousis, Sofia Kokonezi, and John Livieratos.
\newblock Algorithmically efficient syntactic characterization of possibility
  domains.
\newblock {\em Bulletin of the Hellenic Mathematical Society (HMS)},
  63:97--135, 2019.

\bibitem{dokow2009aggregation}
Elad Dokow and Ron Holzman.
\newblock Aggregation of binary evaluations for truth-functional agendas.
\newblock {\em Social Choice and Welfare}, 32(2):221--241, 2009.

\bibitem{dokow2010aggregation}
Elad Dokow and Ron Holzman.
\newblock Aggregation of binary evaluations.
\newblock {\em Journal of Economic Theory}, 145(2):495--511, 2010.

\bibitem{dokow2010aggregationnonB}
Elad Dokow and Ron Holzman.
\newblock Aggregation of non-binary evaluations.
\newblock {\em Advances in Applied Mathematics}, 45(4):487--504, 2010.

\bibitem{enderton2001mathematical}
Herbert~B. Enderton.
\newblock {\em A mathematical introduction to logic}.
\newblock Elsevier, 2001.

\bibitem{DBLP:reference/choice/Endriss16}
Ulle Endriss.
\newblock Judgment aggregation.
\newblock In Felix Brandt, Vincent Conitzer, Ulle Endriss, J{\'{e}}r{\^{o}}me
  Lang, and Ariel~D. Procaccia, editors, {\em Handbook of Computational Social
  Choice}, pages 399--426. Cambridge University Press, 2016.

\bibitem{endriss2016succinctness}
Ulle Endriss, Umberto Grandi, Ronald De~Haan, and J{\'e}r{\^o}me Lang.
\newblock Succinctness of languages for judgment aggregation.
\newblock In {\em Fifteenth International Conference on the Principles of
  Knowledge Representation and Reasoning}, 2016.

\bibitem{endriss2012complexity}
Ulle Endriss, Umberto Grandi, and Daniele Porello.
\newblock Complexity of judgment aggregation.
\newblock {\em Journal of Artificial Intelligence Research}, pages 481--514,
  2012.

\bibitem{fishburn1986aggregation}
Peter~C Fishburn and Ariel Rubinstein.
\newblock Aggregation of equivalence relations.
\newblock {\em Journal of classification}, 3(1):61--65, 1986.

\bibitem{grandi2013lifting}
Umberto Grandi and Ulle Endriss.
\newblock Lifting integrity constraints in binary aggregation.
\newblock {\em Artificial Intelligence}, 199:45--66, 2013.

\bibitem{DBLP:books/daglib/0095988}
Neil Immerman.
\newblock {\em Descriptive complexity}.
\newblock Graduate texts in computer science. Springer, 1999.

\bibitem{jeavons1995algebraic}
Peter Jeavons and David Cohen.
\newblock An algebraic characterization of tractable constraints.
\newblock In {\em International Computing and Combinatorics Conference}, pages
  633--642. Springer, 1995.

\bibitem{jeavons1999determine}
Peter Jeavons, David Cohen, and Marc Gyssens.
\newblock How to determine the expressive power of constraints.
\newblock {\em Constraints}, 4(2):113--131, 1999.

\bibitem{kasher1997question}
Asa Kasher and Ariel Rubinstein.
\newblock On the question ``who is a {J}?" a social choice approach.
\newblock {\em Logique et Analyse}, pages 385--395, 1997.

\bibitem{DBLP:conf/RelMiCS/KirousisKL18}
Lefteris Kirousis, Phokion~G. Kolaitis, and John Livieratos.
\newblock On the computational complexity of non-dictatorial aggregation.
\newblock In {\em Relational and Algebraic Methods in Computer Science - 17th
  International Conference (RAMiCS)}, pages 350--365, 2018.

\bibitem{kirousis2019aggregation}
Lefteris Kirousis, Phokion~G Kolaitis, and John Livieratos.
\newblock Aggregation of votes with multiple positions on each issue.
\newblock {\em ACM Transactions on Economics and Computation (TEAC)}, 7(1):1,
  2019.

\bibitem{lang2002conditional}
J{\'e}r{\^o}me Lang, Paolo Liberatore, and Pierre Marquis.
\newblock Conditional independence in propositional logic.
\newblock {\em Artificial Intelligence}, 141(1-2):79--121, 2002.

\bibitem{lang2003propositional}
J{\'e}r{\^o}me Lang, Paolo Liberatore, and Pierre Marquis.
\newblock Propositional independence-formula-variable independence and
  forgetting.
\newblock {\em Journal of Artificial Intelligence Research}, 18:391--443, 2003.

\bibitem{lang2016agenda}
J{\'e}r{\^o}me Lang, Marija Slavkovik, and Srdjan Vesic.
\newblock Agenda separability in judgment aggregation.
\newblock In {\em Thirtieth AAAI Conference on Artificial Intelligence}, 2016.

\bibitem{larose2017algebra}
Beno{\^\i}t Larose.
\newblock Algebra and the complexity of digraph {CSP}s: a survey.
\newblock In {\em Dagstuhl Follow-Ups}, volume~7. Schloss
  Dagstuhl-Leibniz-Zentrum fuer Informatik, 2017.

\bibitem{DBLP:books/sp/Libkin04}
Leonid Libkin.
\newblock {\em Elements of Finite Model Theory}.
\newblock Texts in Theoretical Computer Science. An {EATCS} Series. Springer,
  2004.

\bibitem{list2002aggregating}
Christian List and Philip Pettit.
\newblock Aggregating sets of judgments: An impossibility result.
\newblock {\em Economics \& Philosophy}, 18(1):89--110, 2002.

\bibitem{list2009judgment}
Christian List and Clemens Puppe.
\newblock Judgment aggregation: A survey.
\newblock {\em Handbook of Rational and Social Choice}, pages 457--482, 2009.

\bibitem{mongin2008factoring}
Philippe Mongin.
\newblock Factoring out the impossibility of logical aggregation.
\newblock {\em Journal of Economic Theory}, 141(1):100--113, 2008.

\bibitem{nehring2002stategy}
Klaus Nehring and Clemens Puppe.
\newblock Strategy-proof social choice on single-peaked domains: Possibility,
  impossibility and the space between, 2002.
\newblock University of California at Davis; available at:
  http://vwl1.ets.kit.edu/puppe.php.

\bibitem{nehring2007structure}
Klaus Nehring and Clemens Puppe.
\newblock The structure of strategy-proof social choice—part i: General
  characterization and possibility results on median spaces.
\newblock {\em Journal of Economic Theory}, 135(1):269--305, 2007.

\bibitem{nehring2010abstract}
Klaus Nehring and Clemens Puppe.
\newblock Abstract arrowian aggregation.
\newblock {\em Journal of Economic Theory}, 145(2):467--494, 2010.

\bibitem{post1941two}
Emil~Leon Post.
\newblock {\em The two-valued iterative systems of mathematical logic}.
\newblock Number~5 in Annals of Mathematics Studies. Princeton University
  Press, 1941.

\bibitem{rubinstein1986algebraic}
Ariel Rubinstein and Peter~C Fishburn.
\newblock Algebraic aggregation theory.
\newblock {\em Journal of Economic Theory}, 38(1):63--77, 1986.

\bibitem{schaefer1978complexity}
Thomas~J. Schaefer.
\newblock The complexity of satisfiability problems.
\newblock In {\em Proc.\ of the 10th Annual ACM Symp.\ on Theory of Computing},
  pages 216--226, 1978.

\bibitem{sharir1981strong}
Micha Sharir.
\newblock A strong-connectivity algorithm and its applications in data flow
  analysis.
\newblock {\em Computers \& Mathematics with Applications}, 7(1):67--72, 1981.

\bibitem{stockmeyer1976polynomial}
Larry~J Stockmeyer.
\newblock The polynomial-time hierarchy.
\newblock {\em Theoretical Computer Science}, 3(1):1--22, 1976.

\bibitem{szegedy2015impossibility}
Mario Szegedy and Yixin Xu.
\newblock Impossibility theorems and the universal algebraic toolkit.
\newblock {\em CoRR}, abs/1506.01315, 2015.

\bibitem{szendrei1986clones}
{\'A}gnes Szendrei.
\newblock {\em Clones in universal algebra}, volume~99.
\newblock Presses de l'Universit{\'e} de Montr{\'e}al, 1986.

\bibitem{Tarjan72}
Robert~Endre Tarjan.
\newblock Depth-first search and linear graph algorithms.
\newblock {\em {SIAM} J. Comput.}, 1(2):146--160, 1972.

\bibitem{terzopoulou2018aggregating}
Zoi Terzopoulou, Ulle Endriss, and Ronald de~Haan.
\newblock Aggregating incomplete judgments: Axiomatisations for scoring rules.
\newblock In {\em Proceedings of the 7th international workshop on
  computational social choice (COMSOC)}, 2018.

\bibitem{zanuttini2002unified}
Bruno Zanuttini and Jean-Jacques H{\'e}brard.
\newblock A unified framework for structure identification.
\newblock {\em Information Processing Letters}, 81(6):335--339, 2002.

\end{thebibliography}

\end{document}